\documentclass{article}
\usepackage{amsfonts,mathenv,amssymb,amsmath}
\usepackage{color}
\numberwithin{equation}{section}
\newtheorem{lemma}{Lemma}[section]
\newtheorem{definition}[lemma]{Definition}
\newtheorem{proposition}[lemma]{Proposition}
\newtheorem{corollary}[lemma]{Corollary}
\newtheorem{theorem}[lemma]{Theorem}
\DeclareMathOperator{\divg}{div}
\DeclareMathOperator{\curl}{curl}
\DeclareMathOperator{\vol}{vol}
\newcommand{\calt}{\mathcal T}

\newcommand{\fra}{\mathfrak a}
\title{Free boundary problem of low Mach number magnetohydrodynamic flows\thanks{ Submitted to the editors DATE.}}
\author{Ruixi Zhang\thanks{Department of Mathematical Sciences, Tsinghua University, Beijing 100084, China (zhangrx24@mails.tsinghua.edu.cn).}}
\date{}
\begin{document}
\maketitle

\begin{abstract}
In this paper we consider a free boundary problem of low Mach number magnetohydrodynamic flow in spatial dimension $n\geq2$. A priori estimates of the second fundamental form and various flow quantities in Sobolev norms are obtained by adopting the geometrical point of view introduced by Christodoulou and Lindblad \cite{CLin}. Moreover, a blow up criterion is derived by using the method of Beale, Kato and Majda \cite{BKM}. 
\end{abstract}

\section{Introduction}\label{sec:intro}
We consider the motion of low Mach number magnetohydrodynamic (MHD) flows in vacuum \cite{JJLX, LZeng}: 
\begin{equation}\label{eq:LM}
\begin{aligned}
&\calt^{-1}(\partial_t+v^k\partial_k)v-B^k\partial_kB=-\partial q\qquad\mathrm{in}\ \mathcal D,\\
&(\partial_t+v^k\partial_k)B-B^k\partial_kv=-(\divg v)B,\qquad\mathrm{in }\ \mathcal D,\\
&(\partial_t+v^k\partial_k)\calt=\calt\Delta\calt\qquad\mathrm{in }\ \mathcal D\\
&\divg v=\Delta\calt,\quad\divg B=0\qquad\mathrm{in }\ \mathcal D. 
\end{aligned}
\end{equation}
Here $\mathcal D\subset[0,T]\times\mathbb R^n,\ n\geq2$ is the domain occupied by the fluid, $v=(v^1,...,v^n)$ is the fluid velocity, $B=(B^1,...,B^n)$ is magnetic field, $q=p+\frac{1}{2}|B|^2$ is the total pressure (where $p$ is the fluid pressure), and $\calt$ is the temperature. Throughout this paper, we use $D_t=\partial_t+v^k\partial_k$ to denote the material derivative. 
The fluid should satisfy the following boundary conditions for $(t,x)\in\partial\mathcal D$:  
\begin{align}
&D_t\in \mathbf{T}\partial\mathcal D,\label{eq:Dttan}\\
&B\cdot N=0,\label{eq:pfcond}\\
&\calt=1\label{eq:diriT},\\
&q=0\text{ and }-\partial_N q=-N^j\partial_jq>0,\label{eq:nosurftens}
\end{align}
where $\mathbf{T}\partial\mathcal D$ is the tangential bundle of $\partial\mathcal D$ and $N$ is the outer normal to the boundary of $\mathcal D\cap\{t\}\times\mathbb R^n$. Here (\ref{eq:Dttan}) means that the free boundary moves with the fluid velocity; (\ref{eq:pfcond}) means that the fluid is a perfect conductor. 
In (\ref{eq:nosurftens}), $q=0$ means that there is no surface tension and the requirement $-\partial_N q>0$ is a natural extension of the Taylor sign condition into MHD, in absence of which free boundary problems for MHD can be ill-posed \cite{Ebin, HaoL}. 

Given a simply connected $\mathcal D_0\subset\mathbb R^n$ and a triple $(v_0,B_0,\calt_0)$ satisfying $\divg B_0=0$, $\divg v_0=\Delta\calt_0$ and (\ref{eq:pfcond}), (\ref{eq:diriT}) on $\partial\mathcal D_0$, we want to find a set $\mathcal D\subset[0,T]\times\mathbb R^n$ and functions $(v,B,\calt)$ solving (\ref{eq:LM}), subject to boundary conditions (\ref{eq:Dttan})-(\ref{eq:nosurftens}) and the initial condition: 
\begin{equation}\label{eq:inicond}
\mathcal D\cap\{t=0\}=\mathcal D_0,\quad(v,\calt,B)|_{t=0}=(v_0,\calt_0,B_0),
\end{equation}
For simplicity, we denote $\mathcal D_t=\mathcal D\cap\{t\}\times\mathbb R^n$. Before moving into a survey of related research, we leave several remarks. Firstly, once $\divg B_0=0$ and $B_0\cdot N=0$ are satisfied by the initial data, they automatically propagate to any positive time, due to the second equation in (\ref{eq:LM}) together with (\ref{eq:Dttan}). Secondly, the function $q$ is determined from $(v,B,\calt)$ by the boundary condition $q|_{\partial\mathcal D_t}=0$ and the elliptic equation obtained from taking divergence on both sides of the first equation in (\ref{eq:LM}). Finally, our boundary condition (\ref{eq:diriT}) is a preassumed temperature distribution on the free surface. It can be replaced by conditions on the heat flux, $\partial_N\calt=0$ \cite{GWang2}, or the heat transfer coefficient, $k_1\partial_N\calt=k_2(\calt_b-\calt)$ \cite{HLuo}.    

The system (\ref{eq:LM}) is derived by letting the Mach number tend to zero in the heat-conductive MHD equation; for a derivation in a viscous setting, see \cite{JJLX}. In that paper, the authors establish uniform-in-Mach-number estimates of solutions to the full compressible MHD equation on $\mathbb R^3$, which is used to prove convergence of solutions as the Mach number goes to zero. 
Now let us mention some remarkable works on free surface motion of fluids with close relation to our present study. 
For the incompressible Euler equation, well-posedness theories are first obtained for irrotational water waves; see \cite{Wu1,Wu2} for local existence theory and \cite{Wu3,Wu4} for global behavior of solutions. 
When vorticity is taken into account, Christodoulou and Lindblad \cite{CLin} prove a priori estimates for the incompressible Euler equation in Sobolev spaces by introducing an energy structure which is able to control the second fundamental form of the free surface. Later, in \cite{Lind,CShk,ZZh} local well-posedness for the incompressible Euler equation is studied. For well-posedness at low regularity and breakdown criterion, see \cite{Gins,TatEul,WZh21}. 
For the compressible inviscid flow, well-posedness theories are first studied for liquids in \cite{Lind2,Trak}. As for gases in a physical vacuum, local well-posedness theory is established in \cite{CShk2,JMas} and global behavior of solutions is studied in \cite{HJang,Zeng}. 

When the flow is heat-conducting and highly subsonic, its motion is described by the low Mach number system, which has appealed much research interest in recent years. Luo and Zeng \cite{LZeng} prove a priori estimate in Sobolev spaces for the low Mach number Euler-Fourier system under the boundary condition (\ref{eq:diriT}), by extending the energy structure in \cite{CLin}. Later, Gu and Wang \cite{GWang2} prove local well-posedness and uniform-in-Mach-number estimate for the Euler-Fourier system under the boundary condition $\partial_N\calt=0$ by extending the method formulated in \cite{CShk}. This allows them to justify convergence towards the low Mach number limit system. 

The magnetohydrodynamic equation describes the motion of conducting fluids in an electromagnetic field, and the corresponding free surface problem is a special case of the classical plasma-vacuum interface problem \cite{GPoe}. The strong coupling between the velocity and magnetic fields leads to challenges in mathematical analysis. 
In \cite{HLuo}, a priori estimate for the incompressible MHD equation is proven under a more restrictive boundary condition, $p=0,\ |B|=\mathrm{const}$ on $\partial\mathcal D_t$, in comparison with the first equation in (\ref{eq:nosurftens}); see also \cite{FHYZh} for a continuation criterion for the problem considered in \cite{HLuo}. 
On the well-posedness theories, local-in-time solutions are constructed in \cite{GWang} for incompressible flow and in \cite{TWang} for compressible flow, both under the zero surface tension condition (\ref{eq:nosurftens}). 
More recently, for the incompressible MHD equation, local well-posedness at sharp regularity is established under the zero surface tension condition in \cite{IPTT}, and local existence and uniqueness of solution is proven under the nonzero surface tension condition in \cite{GLZh}. 
 
In contrast with the extensive study on the MHD system and the low Mach number Euler-Fourier system, there seems to be a research gap in the free surface problem of the coupled system (\ref{eq:LM}). 
In this paper, we shall establish a priori energy estimate for (\ref{eq:LM}), using the method introduced in \cite{CLin}. Moreover, we develop a criterion for the solution to blow up in Sobolev spaces, using the method in \cite{BKM}. 

Before a precise statement of our main results, let us briefly mention what happens when the magnetic field is added into the low Mach number system considered in \cite{LZeng}. The main difference is that, in a word, the transport structure splits into a symmetric hyperbolic structure. For example, in (\ref{eq:LM}) $D_tv$ is not directly related to $\partial q$, and $\curl v$ is no longer transported along $D_t$. Instead, the first two equations in (\ref{eq:LM}) should be coupled together for cancellation of the ``convective terms". To deal with this effect, apart from combining the energy structures in \cite{HLuo, LZeng}, we also have to split the functional $\int_{\mathcal D_t}|\partial D_t^{r-1}\divg v|^2dx$ introduced in \cite{LZeng} into two pieces. Indeed, this functional corresponds to the equation 
\[
D_t\divg v-\calt\Delta\divg v=[D_t,\Delta]\calt+2\nabla\calt\cdot\nabla\divg v+(\divg v)^2, 
\]
which is discovered therein and still holds true for our (\ref{eq:LM}). After applying $D_t$, the leading term of the right-hand side is
\[
D_t[D_t,\Delta]\calt=D_t\curl\curl v^k\partial_k\calt+\cdots,
\]
according to the decomposition $\Delta v=\partial\divg v-\curl\curl v$ (holds for $n=3$). Hence, we have to estimate $D_t\curl v$, which is fine when $B=0$. However, when $B\neq0$, in view of (\ref{eq:curlv}) this involves the ``convective term" $B^i\partial_i\curl B$, which is unbounded. For this problem, our strategy is to match the characteristics $D_\mp$ with the Riemann invariants $v\pm\calt^{1/2}B$. In fact, $D_\mp(v\pm\calt^{1/2}B)$ 
enjoys better estimates than $D_tv,D_tB$, especially after applying $D_t$ and the $\curl$ operator; see sections 3.1, 6.1. 
This motivates our construction of energy (\ref{eq:Er}), as stated below. 

\subsection{Main results}\label{sub:mainres}

The construction of energy involves a bilinear form on $(0,r)$-tensors whose restriction on $\partial\mathcal D_t$ is induced by the tangential projection. To define this bilinear form, we introduce the following
\begin{definition}
The injectivity radius $\iota_0$ of the normal exponential map of $\partial\mathcal D_t$ is the largest number such that the map
\[
\partial\mathcal D_t\times(-\iota_0,\iota_0)\to\mathbb R^n:\ (\bar x,s)\mapsto \bar x+sN(\bar x)
\]
is injective. 
\end{definition}

\begin{definition}
Let $d_0=\frac{1}{16}\iota_0(\mathcal D_0)$ be a fixed number and $\eta$ be a smooth cutoff function which equals 1 on $(-\infty,d_0]$ and vanishes on $[2d_0,\infty)$. Let the outer normal $N$ to $\partial\mathcal D_t$ be extended to the interior by
\begin{equation}\label{eq:Ndef}
N^i(t,x)=\delta^{ij}N_i(t,x),\quad N_i(t,x)=-\eta(d(t,x))\partial_id(t,x),
\end{equation}
where $d(t,x)=d(x,\partial\mathcal D_t)$. 
\end{definition}
With this setting, the bilinear form is specified as
\begin{equation}\label{eq:Qdef}
Q(\alpha,\beta)=Q^{i_1j_1}\cdots Q^{i_rj_r}\alpha_{i_1\cdots i_r}\beta_{j_1\cdots j_r},\qquad Q^{ij}=\delta^{ij}-N^iN^j. 
\end{equation}
Now we are able to construct the energy functional $E_r(t)$ for 
\[
r>r_*=[\frac{n}{2}]+2.
\]
Define
\begin{equation}\label{eq:Er}
E_r(t)=E_r^a(t)+E_r^b(t)+E_r^{c+}(t)+E_r^{c-}(t),
\end{equation}
where
\begin{equation}\label{eq:Era}
\begin{aligned}
E_r^a(t)=\sum_{k=1}^n\int_{\mathcal D_t}\calt^{-1}Q(\partial^rv_k,\partial^rv_k)dx+\sum_{k=1}^n\int_{\mathcal D_t}Q(\partial^rB_k,\partial^rB_k)dx\\
+\int_{\mathcal \partial D_t}Q(\partial^rq,\partial^rq)\mathfrak adS,
\end{aligned}
\end{equation}
\begin{equation}\label{eq:Erb}
E_r^b(t)=\int_{\mathcal D_t}\calt^{-1}|\partial^{r-1}\curl v|^2+|\partial^{r-1}\curl B|^2dx,
\end{equation}
\begin{equation}\label{eq:Erc}
\begin{aligned}
&E_r^{c+}(t)=\int_{\mathcal D_t}|\partial D_t^{r-2}D_-\divg(v+2b)|^2dx,\\
&E_r^{c-}(t)=\int_{\mathcal D_t}|\partial D_t^{r-2}D_+\divg(v-2b)|^2dx, 
\end{aligned}
\end{equation}
and
\[
\mathfrak a=\frac{1}{-\partial_N q},\qquad b=\calt^\frac{1}{2}B,\qquad D_\pm=D_t\pm b^k\partial_k. 
\]
Also define
\begin{equation}\label{eq:calEr}
\mathcal E_r(t)=\mathcal E_{r_*}(t)+\sum_{k=r_*+1}^rE_k(t),
\end{equation}
where $\mathcal E_{r_*}$ is the square integral of derivatives of $(v,B,q,\calt,\divg v)$ up to order $r_*$; see (\ref{eq:Er*def}). With these notations, our first result is stated as follows: 

\begin{theorem}\label{thm:main}
Suppose that $(\mathcal D_t,v,B,q,\calt)$ is a smooth solution to (\ref{eq:LM}) and that the initial data satisfies
\begin{align}
-\partial_Nq&\geq\epsilon_{0*}\qquad\text{on }\partial\mathcal D_0,\label{eq:TS0}\\
\calt&\geq\epsilon_{1*}\qquad\text{in }\mathcal D_0,\label{eq:infT0}\\
|\theta|+\frac{1}{\iota_0}&\leq K_*\qquad\text{on }\partial\mathcal D_0,\label{eq:iniass1}\\
|\partial^2q|+|\partial D_tq|&\leq L_*\qquad\text{on }\partial\mathcal D_0,\label{eq:iniass2}
\end{align}
Then there is a $T=T(\epsilon_{0*},\epsilon_{1*},K_*,L_*,\vol\mathcal D_0,\mathcal E_{r}(0))>0$ such that the following estimates hold for $t\in[0,T]$:
\begin{align}
-\partial_Nq&\geq\frac{1}{2}\epsilon_{0*}\qquad\text{on }\partial\mathcal D_t,\label{eq:recTS}\\
\calt&\geq\epsilon_{1*}\qquad\text{in }\mathcal D_t,\label{eq:recinfT}\\
|\theta|+\frac{1}{\iota_0}&\leq 2K_*\qquad\text{on }\partial\mathcal D_t,\label{eq:recass1}\\
|\partial^2q|+|\partial D_tq|&\leq 2L_*\qquad\text{on }\partial\mathcal D_t,\label{eq:recass2}
\end{align}
\begin{equation}\label{eq:volest}
\frac{1}{2}\vol\mathcal D_0\leq\vol\mathcal D_t\leq 2\vol\mathcal D_0
\end{equation}
and
\begin{equation}\label{eq:Eest}
\mathcal E_r(t)\leq 2\mathcal E_r(0).
\end{equation} 
\end{theorem}
%The theorem states that the energy functional $\mathcal E_r$ is stable under the evolution of (\ref{eq:LM}). In turn, $\mathcal E_r$ is able to control the regularity of the solution. 
As a result, the regularity of the second fundamental form and the flow variables can be controlled, as stated below: 
\begin{corollary}\label{cor:regstab}
Let $(\mathcal D_t,v,B,q,\calt)$ and $T$ be as above. Then we have
\[\begin{aligned}
\sum_{l=0}^{r-2}\sum_{k=0}^{r-l}\Vert\partial^k D_t^l(v,B,q,\calt,\divg v)\Vert_{L^2(\mathcal D_t)}^2&+\sum_{k=0}^{r-2}\Vert\overline\partial^k\theta\Vert_{L^2(\partial\mathcal D_t)}^2\\
&\leq C(\epsilon_{0*},\epsilon_{1*},K_*,L_*,\vol\mathcal D_0,\mathcal E_{r}(0)). 
\end{aligned}\]
\end{corollary}
\textbf{Remark. }In Theorem \ref{thm:main} we do not impose initial bounds of $|\partial v|, |\partial q|$ and $|D_t\divg v|$, in contrast to \cite{CLin} and \cite{LZeng}, because they are already controlled by $\mathcal E_{r_*}(0)$. \\[2mm]
\textbf{Remark. }In \cite{LZeng}, the recovered bound of $|\theta|+\frac{1}{\iota_0}$ is $CK_*$ for some certain constant $C$, in contrast to (\ref{eq:recass2}) above. We improve this bound by quantifying the loss of $\iota_0$ in terms of $|D_t\theta|$ and the Lipschitz distance between $\partial\mathcal D_0$ and $\partial\mathcal D_t$. \\[2mm]
Our second result is the following blow up criterion in dimension $n=3$: 
\begin{theorem}\label{thm:buc}
Suppose $n=3$ and $(\mathcal D_t,v,B,q,\calt)$ is a smooth solution to (\ref{eq:LM}) for $t\in[0,T)$. Define
\begin{equation}\label{eq:Mdef}
M=\Vert\divg v\Vert_{L^\infty(\mathcal D_t)}+\Vert\curl v\Vert_{L^\infty(\mathcal D_t)}+\Vert\curl B\Vert_{L^\infty(\mathcal D_t)}+\Vert\overline\partial(v,B)\Vert_{L^\infty(\partial\mathcal D_t)},
\end{equation}
\begin{equation}\label{eq:Kdef}
K=\Vert\theta\Vert_{L^\infty(\partial\mathcal D_t)}+\frac{1}{\iota_0(\partial\mathcal D_t)}. 
\end{equation}
Suppose 
\[
-\partial_N q\geq\epsilon_0>0\quad\text{on }\partial\mathcal D_t,
\]
and suppose
\begin{equation}\label{eq:MKbar}
\limsup_{t\to T}M=\overline M<\infty,\qquad\limsup_{t\to T}K=\overline K<\infty,
\end{equation}
\begin{equation}\label{eq:DNDtq}
\int_0^T\Vert\partial_N D_tq\Vert_{L^\infty(\partial\mathcal D_t)}dt=\overline L<\infty,
\end{equation}
then the solution is uniformly bounded in the sense that
\[
\sup_{t\in[0,T)}\mathcal E_r(t)\leq C(T,\epsilon_0,\overline K,\overline M,\overline L,L_0,L_1,\mathcal E(0))\mathcal E_r(0),\quad r\geq4,
\]
where $L_0,L_1$ only depend on the initial data, as defined in (\ref{eq:L0def}), (\ref{eq:L1def}), and $\mathcal E$ is defined in (\ref{eq:calEdef}). Moreover, when $B=0$ the condition (\ref{eq:DNDtq}) can be removed. 
\end{theorem}
\textbf{Remark. }The energy estimate in \cite{LZeng} does not imply Theorem \ref{thm:buc} for $B=0$. Indeed, in \cite{LZeng} a priori assumptions are imposed on the $L^\infty(\mathcal D_t)$ bound of $\partial v$, $\partial\divg v$ and $D_t\divg v$, which cannot be controlled by $M$ in (\ref{eq:Mdef}). Our estimate is more delicate in several aspects. Firstly, we adopt the method in \cite{BKM} for estimating $\Vert\partial v\Vert_{L^\infty(\mathcal D_t)}$. Secondly, we frequently use interpolation inequalities for nonlinear estimates. Thirdly, we derive a parabolic equation of $D_t\divg v$ in the form (\ref{eq:pardiv}). This enables us to obtain a priori estimates of $D_t\divg v$, $\partial\divg v$. \\[2mm]
\textbf{Remark. }Theorem \ref{thm:buc} roughly reduces to the criterion in \cite{Gins} if $\calt$ is a constant and $B=0$, except that we need $\limsup M(t)$ instead of $\int M(t)^2dt$. Moreover, there is no Dirichlet-to-Neumann operator in the definition of $M$. This improvement relies on the gradient estimate for harmonic functions. The high nonlinearity in $M$ is mainly due to the presence of $\partial\calt$, $\divg v$. \\[2mm]
Let us mention some new difficulties in deriving Theorem \ref{thm:buc}. 
Compared with \cite{LZeng}, the presence of $B$ forces us to match the characteristics $D_\mp$ with the Riemann invariants $v\pm\calt^{1/2}B$, as discussed before. This is important not only for recovering one more derivative, but also for controlling the nonlinearity in the energy estimate. 
The proofs in \cite{FHYZh,Gins} are based on a nonlinear estimate of the form 
\begin{equation}\label{eq:diffineq}
\frac{d}{dt}\mathcal E(t)\lesssim \mathcal E(t)(1+\log^+\mathcal E(t)). 
\end{equation}
This combined with $\int_1^\infty[y(1+\log y)]^{-1}dy=\infty$ prohibits the blowup of $\mathcal E(t)$. The proof of Theorem \ref{thm:buc} is in a similar spirit. However, when $B\neq0$ the estimate of $\frac{d}{dt}\mathcal E(t)$ involves a $(\log\mathcal E)^2$ factor, which is unacceptable. To overcome this problem, we replace (\ref{eq:diffineq}) by its integral analogue: 
\begin{equation}\label{eq:intineq}
\mathcal E(t)\lesssim \mathcal E(0)+\int_0^t\mathcal E(s)(1+\log^+\mathcal E(s))ds. 
\end{equation}
Using integration by parts in spacetime, we are able to lower down the nonlinearity in $\mathcal E$ and obtain (\ref{eq:intineq}).

\section{Preliminaries}

\subsection{Geometric aspects}

In \cite{CLin}, a Lagrangian formulation is used. Denote $\Omega=\mathcal D_0$ and define a coordinate change $\Omega\ni(t,y)\mapsto(t,x)\in\mathcal D_t$ by 
\begin{equation}\label{eq:lagcoor}
\frac{d}{dt}x(t,y)=v(t,x(t,y)),\quad x(0,y)=y,
\end{equation}
so that the temporal derivetive $\partial_t|_{y=\mathrm{const}}$ in $(t,y)$ coordinate coincides with the material deivative $D_t$ defined in section \ref{sec:intro}. 
For fixed $t$, the Euclidean metric $\delta_{ij}$ pulled back by the map $y\mapsto x(t,y)$ leads to a metric $g_{ab}(t,y)$ on $\Omega$ by
\[
g_{ab}(t,y)=\delta_{ij}\frac{\partial x^i}{\partial y^a}\frac{\partial x^j}{\partial y^b}. 
\]
Tensors are pulled back by 
\begin{equation}\label{eq:tensorlaw}
\alpha_{a_1\cdots a_r}=\frac{\partial x^{i_1}}{\partial y^{a_1}}\cdots\frac{\partial x^{i_r}}{\partial y^{a_r}}\alpha_{i_1\cdots i_r}, 
\end{equation}
and their length is defined by
\begin{equation}\label{eq:glength}
|\alpha_{a_1\cdots a_r}|=\left(g^{a_1b_1}\cdots g^{a_rb_r}\alpha_{a_1\cdots a_r}\alpha_{b_1\cdots b_r}\right)^\frac{1}{2}, 
\end{equation}
where $g^{ab}$ is the inverse matrix of $g_{ab}$, that is, $g^{ab}=\delta^{ij}\frac{\partial y^a}{\partial x^i}\frac{\partial y^b}{\partial x^j}$. %The volume measure corresponding to $g$ is $\sqrt{\det(g_{ab})}dy$. 
Now that $t\mapsto x(t,y)$ is isometric, it is convenient to identify the geometries in $\Omega$ and $\mathcal D_t$. For example, if we denote by $\nabla_i$ the covariant derivative in the direction $\frac{\partial}{\partial x^i}$, then
\begin{equation}\label{eq:covder}
\nabla_if(x(t,y))=\frac{\partial f}{\partial x^i}(x(t,y)),\qquad\nabla_i\alpha_{i_1\cdots i_r}=\frac{\partial\alpha_{i_1\cdots i_r}}{\partial x^i}. 
\end{equation}
From this, we introduce the convention
\begin{equation}\label{eq:Eulind}
\textit{Hereinafter, tensors are expressed with respect to Eulerian frame in indices }i,j,k...
\end{equation}
and with respect to Lagrangian frame in indices $a,b,c...$. 

Let us move on to look at the boundary geometry. Recall the extended outer normal $N$ from (\ref{eq:Ndef}). Define the tangential projection operator
\begin{equation}\label{eq:Pidef}
\begin{aligned}
\Pi(\alpha)_{i_1\cdots i_r}=&\Pi_{i_1\cdots i_r}^{j_1\cdots j_r}\alpha_{j_1\cdots j_r}\\
=&(\delta_{i_1}^{j_1}-N_{i_1}N^{j_1})\cdots(\delta_{i_r}^{j_r}-N_{i_r}N^{j_r})\alpha_{j_1\cdots j_r} 
\end{aligned}
\end{equation}
and note the identity
\begin{equation}\label{eq:QPi}
Q(\alpha,\beta)=\langle\Pi(\alpha),\Pi(\beta)\rangle=\delta^{IJ}\Pi(\alpha)_I\Pi(\beta)_J. 
\end{equation}
We summarize some useful results below; see Lemmas 2.1, 3.9, 3.10 in \cite{CLin}. 

\begin{lemma}\label{lem:Nmu}
(i) Let $\mu_\zeta$ be the area measure induced by $\zeta_{ij}=\delta_{ij}-N_iN_j$ and $\mu_g$ be the volume measure induced by $g$. Then we have
\begin{equation}\label{eq:Dtmu}
D_td\mu_g=\divg v\,d\mu_g,\qquad D_td\mu_\zeta=(\divg v-N^iN^j\nabla_iv_j)\,d\mu_\zeta. 
\end{equation}
(ii) Suppose $|\theta|+\frac{1}{\iota}\leq K$. Then we have
\begin{equation}\label{eq:gradN}
\Vert\nabla N\Vert_{L^\infty(\Omega)}\leq CK,
\end{equation}
\begin{equation}\label{eq:DtNint}
\Vert D_tN\Vert_{L^\infty(\Omega)}\leq C\Vert\nabla v\Vert_{L^\infty(\Omega)}
\end{equation}
and
\begin{equation}\label{eq:DtNbdry}
|D_tN|\leq C|\nabla v|\qquad\text{on }\partial\Omega. 
\end{equation}
\end{lemma}
\textbf{Remark. }In (\ref{eq:Dtmu}) we rewrite the result in \cite{CLin} using Eulerian indices. In (\ref{eq:DtNint}), (\ref{eq:DtNbdry}) one may feel confused about the definition of $D_tN$, because $D_tN_a$ is different from $\frac{\partial x^i}{\partial y^a}D_tN_i$, that is, the pull back (\ref{eq:tensorlaw}) does not commute with $D_t$. However, one readily checks that
\begin{equation}\label{eq:DtNequiv}
\left(g^{ab}D_tN_aD_tN_b\right)^\frac{1}{2}=\left(\delta^{ij}D_tN_iD_tN_j\right)^\frac{1}{2}+O(|\nabla v|),
\end{equation}
thus they make no difference in the estimates (\ref{eq:DtNint}), (\ref{eq:DtNbdry}). \\[2mm]
Define $\overline\nabla_i=\Pi(\nabla)_i=(\delta_i^j-N_iN^j)\nabla_j$. There is a relation between the iterated tangential derivative $\overline\nabla^r$ and the tangential projection of $\nabla^r$: (see Proposition 4.4 in \cite{CLin})

\begin{lemma}
For a scalar function $f$ which vanishes on $\partial\Omega$, it holds on $\partial\Omega$ and for $r\geq2$ that
\begin{equation}\label{eq:proj}
\Pi(\nabla^rf)=(\nabla_Nf)\overline{\nabla}^{r-1}N+\sum_{s=2}^{r-1}\sum_{r_1+\cdots+r_m=r-s}c_{r_1\cdots r_ms}\overline{\nabla}^{r_1}N\cdots\overline{\nabla}^{r_m}N\nabla^sf,
\end{equation}
where $c_{r_1\cdots r_ms}$ is a contraction of $m$ indices. %Moreover, if $\nabla_Nf\nequiv0$, it also holds\begin{equation}\label{eq:Nder}\overline{\nabla}^sN=\sum_{s_1+\cdots+s_m=s}c_{s_1\cdots s_m}(N,|\nabla q|^{-1})\nabla^{s_1+1}q\cdots\nabla^{s_m+1}q,\end{equation}
\end{lemma}
As a corollary, we have
\begin{corollary}\label{cor:DNexpr}
For a scalar function $f$ which vanishes on $\partial\Omega$, it holds on $\partial\Omega$ that
\begin{equation}\label{eq:DNexpr}
\overline\nabla^sN=\sum_{s_1+\cdots+s_m=s}c_{s_1\cdots s_m}(N,|\nabla f|^{-1})\nabla^{s_1+1}f\cdots\nabla^{s_m+1}f,
\end{equation}
where the summation is taken over $s_1,...,s_m\geq1$ and $c_{s_1\cdots s_m}$ is a contraction of $m$ indices whose coefficients depend on $N$ and $|\nabla q|^{-1}$ continuously. 
\end{corollary}
Below we summarize some trace and elliptic estimates; the $L^p$ norms are defined by (\ref{eq:Lpint}), (\ref{eq:Lpbdry}). Inequalities (\ref{eq:bdryHr})-(\ref{eq:diri}) are from Lemma 5.6 in \cite{CLin}, and (\ref{eq:wHr}) is a generalization of (\ref{eq:diri}). 

\begin{lemma}
Let $f$ be any function on $\overline\Omega$ and $r\geq1$. It holds that
\begin{equation}\label{eq:bdryHr}
\Vert\nabla^rf\Vert_{L^2(\partial\Omega)}^2\leq C\left(\int_{\partial\Omega}Q(\nabla^rf,\nabla^rf)d\mu_\zeta+\int_{\Omega}|\nabla^{r-1}\Delta f||\nabla^rf|+K|\nabla^rq|^2d\mu_g\right).
\end{equation}
\end{lemma}
\begin{lemma}%[Lemmas 5.5-5.6 in \cite{CLin}]\label{lem:ellipest}
%Suppose $|\theta|+\iota_0^{-1}\leq K$. 
Let $f$ be any function on $\overline{\Omega}$ and $r\geq2$. We have
\begin{equation}\label{eq:tr}
\Vert f\Vert_{L^2(\partial\Omega)}^2\leq (2\Vert\nabla f\Vert_{L^2(\Omega)}+CK\Vert f\Vert_{L^2(\Omega)})\Vert f\Vert_{L^2(\Omega)},
\end{equation}
\begin{equation}\label{eq:diri}
\begin{aligned}
\Vert\nabla^r f\Vert_{L^2(\Omega)}^2\leq&\Vert N^k\Pi(\nabla^{r-2}\nabla_kf)\Vert_{L^2(\partial\Omega)}\Vert \Pi(\nabla^{r}f)\Vert_{L^2(\partial\Omega)}\\
&+C(\Vert\nabla^{r-2}\Delta f\Vert_{L^2(\Omega)}^2+K^2\Vert\nabla^{r-1}f\Vert_{L^2(\Omega)}^2).
\end{aligned}
\end{equation}
Moreover, let $\rho\geq0$ be a Lipschitz function, then 
\begin{equation}\label{eq:wHr}
\begin{aligned}
\Vert\rho\nabla^rf\Vert_{L^2(\Omega)}^2\leq&\Vert\rho\Vert_{L^\infty(\partial\Omega)}^2\Vert N^k\Pi(\nabla^{r-2}\nabla_kf)\Vert_{L^2(\partial\Omega)}\Vert \Pi(\nabla^{r}f)\Vert_{L^2(\partial\Omega)}\\
&+C(\Vert\rho\nabla^{r-2}\Delta f\Vert_{L^2(\Omega)}^2+K^2\Vert\rho\nabla^{r-1}f\Vert_{L^2(\Omega)}^2+\Vert\nabla\rho\nabla^{r-1}f\Vert_{L^2(\Omega)}^2).
\end{aligned}
\end{equation}
\end{lemma}
\textit{Proof of (\ref{eq:wHr}): }Consider the vector field
\[
Y_k=Q_k^l\delta^{ij}Q(\nabla^{r-2}\nabla_l\nabla_if,\nabla^{r-2}\partial_jf)-Q^{ij}Q(\nabla^{r-2}\nabla_kf,\nabla^{r-2}\nabla_i\nabla_j f),
\]
where $Q_k^l=\delta_{km}Q^{ml}$. Note that $N^kQ_k^l=0$ on $\partial\Omega$, then we have
\[
|N^kY_k|\leq|N^k\Pi(\nabla^{r-2}\nabla_kf)|\cdot|\Pi(\nabla^{r}f)|
\]
and
\[
|\divg Y-\delta^{ij}Q(\nabla^{r-1}\nabla_if,\nabla^{r-1}\nabla_jf)|\leq|\nabla^rf|\big(|\nabla^{r-2}\Delta f|+CK|\nabla^{r-1}f|\big). 
\]
Using
\[
\int_{\partial\Omega}\rho^2Y_kN^kdS=\int_\Omega\rho^2\divg Ydx+\int_\Omega2\rho\nabla\rho\cdot Ydx
\]
and noting
\[
|\nabla^rf|^2\leq \delta^{ij}Q(\nabla^{r-1}\nabla_i f,\nabla^{r-1}\nabla_j f)+C|\nabla^{r-2}\Delta f|^2,
\]
we complete the proof. $\Box$

\subsection{Commutator estimates}

Throughout this paper we specify the $L^p$ norm as
\begin{equation}\label{eq:Lpint}
\Vert f\Vert_{L^p(\Omega)}:=\bigg(\int_\Omega |f|^pd\mu_g\bigg)^{\frac{1}{p}}
\end{equation}
and
\begin{equation}\label{eq:Lpbdry}
\Vert f\Vert_{L^p(\partial\Omega)}:=\bigg(\int_{\partial\Omega} |f|^pd\mu_\zeta\bigg)^{\frac{1}{p}}.
\end{equation}
That is, the measure is always coupled with the metrics $g$ and $\zeta$. The norms $\Vert f\Vert_{W^{k,p}(\Omega)}=\sum_{j=0}^k\Vert\nabla^jf\Vert_{L^p(\Omega)}$ and $\Vert f\Vert_{H^k(\Omega)}=\Vert f\Vert_{W^{k,2}(\Omega)}$ are defined in a usual way. A more frequently used norm is defined as
\begin{equation}\label{eq:rsnorm}
\Vert f\Vert_{r,s}^2:=\sum_{k=0}^r\sum\Vert\Gamma_k\cdots\Gamma_1f\Vert_{L^2(\Omega)}^2,\qquad \Vert f\Vert_{r}:=\Vert f\Vert_{r,0},
\end{equation}
where the summation is taken over $\Gamma_1,...,\Gamma_k\in\{D_t,\nabla\}$ with $D_t$ appearing at most $s$ times. With these notations, we specify $\mathcal E_{r_*}$ in (\ref{eq:calEr}) as
\begin{equation}\label{eq:Er*def}
\mathcal E_{r_*}(t)=\Vert(v,B,\calt,\divg v)\Vert_{r_*,r_*-1}^2+\Vert q\Vert_{r_*,r_*-2}^2+\sum_{k=1}^{r_*}\Vert\overline\nabla^kq\Vert_{L^2(\partial\Omega)}^2. 
\end{equation}
\textbf{Remark. }The norm $\Vert\cdot\Vert_{r,s}$ behaves like $\Vert\cdot\Vert_{H^s(\Omega)}$ under nonlinear actions. For example, we have
\[
\Vert f_1f_2\Vert_{r,s}\leq C(K)\Vert f_1\Vert_{r_1,s_1}\Vert f_1\Vert_{r_1,s_1},\quad r<r_1+r_2-\frac{n}{2},\ s\leq\min\{s_1,s_2\},
\]
as is seen from (\ref{eq:prodint}). \\[2mm]
We present some commutator estimates which will be frequently used. 

\begin{lemma}\label{lem:commest}
Let $S$ and $T$ be compositions of $l$ $D_t$'s and $r-l$ $\nabla$'s. Then
\begin{equation}\label{eq:commform}
S-T=\sum_{m=1}^{r-1}\sum C(\Gamma_1,...,\Gamma_r)(\Gamma_{r}\cdots\Gamma_{m+1}v)\Gamma_m\cdots\Gamma_1,
\end{equation}
where $C(\Gamma_1,...,\Gamma_r)$ is a constant and the summation is taken over $\Gamma_1,...,\Gamma_{r}\in\{D_t,\nabla\}$ with $D_t$ appearing $l-1$ times. As a consequence, if $|\theta|+\frac{1}{\iota_0}\leq K$, then
\begin{equation}\label{eq:[r*]}
\Vert(S-T)f\Vert_{L^2(\Omega)}\leq C(K)\Vert v\Vert_{r_*,l-1}\Vert f\Vert_{r-1,l-1},\quad r\leq r_*,
\end{equation}
\begin{equation}\label{eq:[r]}
\Vert(S-T)f\Vert_{L^2(\Omega)}\leq C(K)\Vert v\Vert_{r-1,l-1}\Vert f\Vert_{r-1,l-1},\quad r>r_*,
\end{equation}
\begin{equation}\label{eq:[r+1]}
\begin{aligned}
\Vert(S-T)f\Vert_{L^2(\Omega)}\leq C(K)(&\Vert v\Vert_{r-1,l-1}\Vert f\Vert_{r-2,l-1}\\
&+\Vert v\Vert_{r-2,l-1}\Vert f\Vert_{r-1,l-1}),\quad r>r_*+1. 
\end{aligned}
\end{equation}
Specially, we have
\begin{equation}\label{eq:rsequiv}
\Vert f\Vert_{r,s}\leq \Vert f\Vert_{r-1,s}+\sum_{l=0}^s\Vert \nabla^{r-l}D_t^lf\Vert_{L^2(\Omega)}+C(K)\Vert v\Vert_{r-1,s-1}\Vert f\Vert_{r-1,s-1},\quad r>r_*.
\end{equation}
\end{lemma}
\textit{Proof: }Suppose $S=\Gamma_1\cdots \Gamma_r$, $\Gamma_i\in\{D_t,\nabla\}$, then there is a permutation $\sigma$ such that $T=\Gamma_{\sigma(1)}\cdots \Gamma_{\sigma(r)}$. Factorize $\sigma$ into interchanges of adjacent indices, then $S-T$ is a sum of the following general terms:
\[
\pm\cdots[D_t,\nabla]\cdots,
\]
where $\cdots$ is a composition of $D_t$ and $\nabla$. By $[\nabla,D_t]=(\nabla v^k)\nabla_k$ and the Leibniz rule, we then obtain (\ref{eq:commform}). Now (\ref{eq:[r*]}), (\ref{eq:[r]}) and (\ref{eq:[r+1]}) are consequences of (\ref{eq:prodint}) and (\ref{eq:rsequiv}) is the corollary of (\ref{eq:[r]}). $\Box$\\[2mm]
A generalization of this lemma is the following
\begin{lemma}\label{lem:commest2}
Let $r\geq s$ and $A_\lambda=a_\lambda^0D_t+a_\lambda^k\nabla_k,\ \lambda=1,...,r$ where $a_1^0=\cdots=a_{r-s}^0=0$. For two permutations $(\lambda_1,...,\lambda_r)$, $(\mu_1,...,\mu_r)$ of $(1,...,r)$ define
\[
S=A_{\lambda_1}\cdots A_{\lambda_r},\quad T=A_{\mu_1}\cdots A_{\mu_r}. 
\]
Suppose $|\theta|+\frac{1}{\iota_0}\leq K$, then we have 
\begin{equation}\label{eq:commr}
\Vert(S-T)f\Vert_{L^2(\Omega)}\leq C(K,\Vert (a,v)\Vert_{r-1,s})\Vert f\Vert_{r-1,s},\quad r>r_*,
\end{equation}
\begin{equation}\label{eq:commr1}
\Vert(S-T)f\Vert_{L^2(\Omega)}\leq C(K,\Vert (a,v,f)\Vert_{r-2,s-1})\Vert (a,v,f)\Vert_{r-1,s},\quad r>r_*+1.
\end{equation}
\end{lemma}
\textit{Proof: }Similar to the proof of Lemma \ref{lem:commest}, we decompose $S-T$ into 
\[
\cdots[A_\lambda, A_\mu]\cdots,
\]
where $\cdots$ is a composition of the $A_\nu$'s. By the calculation
\[
[A_\lambda,A_\mu]=(A_\lambda a_\mu^0-A_\mu a_\lambda^0)D_t+(A_\lambda a_\mu^k-A_\mu a_\lambda^k)\nabla_k+(a_\lambda^0a_\mu^k-a_\lambda^ka_\mu^0)\nabla_k v^l\nabla_l. 
\]
and the Leibniz rule, we further rewrite $(S-T)f$ as a sum of the following general terms: 
\[
C(a,v)\Gamma^{r_1}\phi_1\cdots\Gamma^{r_{m-1}}\phi_{m-1}\Gamma^{r_m}f,\quad r_1+\cdots+r_m=r,\ m\geq2,
\]
where $C(a,v)$ is a polynomial of $a,v$, $\Gamma^{r_i}$ is a composition of $r_i$ elements in $\{D_t,\nabla\}$, $\phi_i\in\{a,v\}$ and $r_1,...,r_m\geq1$. Then (\ref{eq:commr}) follows from (\ref{eq:prodint}). For (\ref{eq:commr1}), we rewrite the above expression as $\Gamma^{r_1}\psi_1\cdots\Gamma^{r_m}\psi_m$ where $\psi_i\in\{a,v,f\}$. Let $j\in\{1,...,m\}$ be such that $\Gamma^{r_j}$ contains most $D_t$'s, then $\Gamma^{r_i}$ contains at most $s-1$ $D_t$'s, $i\neq j$. Due to $r>r_*+2$ we have $\sum_{i\neq j}(r-2-r_i)+(r-1-r_j)>\frac{n}{2}(m-1)$, then (\ref{eq:commr1}) follows from (\ref{eq:prodint}). $\Box$

\begin{lemma}\label{lem:commest3}
Suppose $|\theta|+\frac{1}{\iota_0}\leq K$, then for $k_2+s+k_1=r$ we have
\begin{equation}\label{eq:spcomr}
\Vert\nabla^{k_2}[D_t^s,a\nabla^{k_1}]\Vert_{L^2(\Omega)}\leq C(K,\Vert a\Vert_{r-1,s},\Vert v\Vert_{r-1,s-1})\Vert f\Vert_{r-1,s-1},\quad r>r_*. 
\end{equation}
\end{lemma}
\textit{Proof: }A calculation gives
\[
[D_t^s,a\nabla^{k_1}]f=\sum_{l=1}^s\binom{s}{l}(D_t^la) D_t^{s-l}\nabla^{k_1}f+a[D_t^s,\nabla^{k_1}]f. 
\]
By (\ref{eq:prodint}) we have $\Vert\nabla^{k_2}\{(D_t^la) D_t^{s-l}\nabla^{k_1}f\}\Vert_{L^2(\Omega)}\leq C(K)\Vert a\Vert_{r-1,s}\Vert f\Vert_{r-1,s-1}$. For the last term, we distribute $\nabla^{k_2}$ onto $a,[D_t^s,\nabla^{k_1}]f$ and use (\ref{eq:commform}), (\ref{eq:prodint}). This completes the proof. $\Box$

\section{Proof of Theorem \ref{thm:main}}

We make the following assumptions: 
\begin{align}
-\nabla_Nq&\geq\epsilon_0\qquad\text{on }\partial\Omega,\label{eq:TS}\\
\calt&\geq\epsilon_1\qquad\text{in }\Omega,\label{eq:infT}\\
|\theta|+\frac{1}{\iota_0}&\leq K\qquad\text{on }\partial\Omega,\label{eq:ass1}\\
|\nabla^2q|+|\nabla D_tq|&\leq L\qquad\text{on }\partial\Omega.\label{eq:ass2}
\end{align}
Under these assumptions we shall show that
\begin{equation}\label{eq:dtErest}
\frac{d}{dt}\mathcal E_r(t)\leq C(\mathcal A,\mathcal E_{r-1})\mathcal E_r(t),
\end{equation}
where 
\begin{equation}\label{eq:calAdef}
\mathcal A=(\epsilon_0,\epsilon_1,K,L,\vol\Omega). 
\end{equation}
In sections \ref{sec:dtErc}, \ref{sec:dtErab} we prove (\ref{eq:dtErest}). In section \ref{sec:pfmainfin} we recover the assumptions (\ref{eq:TS})-(\ref{eq:ass2}) and complete the proof of 
Theorem \ref{thm:main}.

\subsection{Estimate of the temporal derivative of $E_r^c$}\label{sec:dtErc}

In this subsection we prove
\begin{equation}\label{eq:dtErcest}
\frac{d}{dt}E_r^{c\pm}(t)\leq C(\mathcal A,\mathcal E_{r-1})\mathcal E_r(t).
\end{equation}
Recall $b=\calt^{\frac{1}{2}}B$ and $D_\pm=D_t\pm b^k\nabla_k$. The first two equations in (\ref{eq:LM}) are rewritten as 
\begin{equation}\label{eq:Dmpvb}
\left\{
\begin{aligned}
D_-(v+b)&=-\calt\nabla q-\frac{1}{2}b\divg(v+2b),\\
D_+(v-b)&=-\calt\nabla q+\frac{1}{2}b\divg(v-2b).
\end{aligned}
\right.
\end{equation}
We need the following parabolic equations for $\divg(v\pm2b)$: 
\begin{lemma}
(i) We have
\begin{equation}\label{eq:Dtdivvb}
\begin{aligned}
D_t\divg(v+2b)-\calt\Delta\divg(v+2b)=[D_+,\Delta]\calt+2\nabla\calt\cdot\nabla\divg(v+2b)+\divg v\,\divg(v+b),\\
D_t\divg(v-2b)-\calt\Delta\divg(v-2b)=[D_-,\Delta]\calt+2\nabla\calt\cdot\nabla\divg(v-2b)+\divg v\,\divg(v-b). 
\end{aligned}
\end{equation}
(ii) For $s\geq1$ we have
\begin{equation}\label{eq:Dtdivvb1}
\begin{aligned}
D_t^sD_-\divg(v+2b)-\calt\Delta(D_t^{s-1}D_-\divg(v+2b))=D_t^{s-1}D_-[D_+,\Delta]\calt+R_s^+,\\
D_t^sD_+\divg(v-2b)-\calt\Delta(D_t^{s-1}D_+\divg(v-2b))=D_t^{s-1}D_+[D_-,\Delta]\calt+R_s^-, 
\end{aligned}
\end{equation}
where
\begin{equation}\label{eq:Rsp}
\begin{aligned}
R_s^+=[D_t^{s-1}D_-,\calt\Delta]\divg(v+2b)+\frac{1}{2}D_t^{s-1}(\divg v\,b^k\nabla_k\divg(v+2b))\\
+D_t^{s-1}D_-(2\nabla\calt\cdot\nabla\divg(v+2b)+\divg v\,\divg(v+b)),
\end{aligned}
\end{equation}
\begin{equation}\label{eq:Rsm}
\begin{aligned}
R_s^-=[D_t^{s-1}D_+,\calt\Delta]\divg(v-2b)-\frac{1}{2}D_t^{s-1}(\divg v\,b^k\nabla_k\divg(v-2b))\\
+D_t^{s-1}D_+(2\nabla\calt\cdot\nabla\divg(v-2b)+\divg v\,\divg(v-b)). 
\end{aligned}
\end{equation}
\end{lemma}
\textit{Proof: }A calculation gives $D_t\calt=\calt\divg v$ and $D_\pm\calt=\calt\divg(v\pm2b)$. Applying $D_\pm$ to $D_t\calt=\calt\Delta\calt$, using
\begin{equation}\label{eq:[b,Dt]}
[b^k\nabla_k,D_t]=\frac{1}{2}\divg v\,b^k\nabla_k
\end{equation}
and dividing both sides by $\calt$, we get
\begin{equation}\label{eq:parmid}
\divg v\divg b+D_t\divg(v+2b)=D_+\Delta\calt,\quad-\divg v\divg b+D_t\divg(v-2b)=D_-\Delta\calt. 
\end{equation}
Commuting derivatives we obtain (\ref{eq:Dtdivvb}). Further applying $D_t^{s-1}D_\mp$ and using (\ref{eq:[b,Dt]}) we obtain (\ref{eq:Dtdivvb1}). $\Box$\\[2mm]
These equations enable us to calculate the temporal derivative of $E_3^{c+}(t)$:  
\[\begin{aligned}
\frac{1}{2}\frac{d}{dt}E_{r}^{c+}(t)=&\frac{1}{2}\int_{\Omega}|\nabla D_t^{r-2}D_-\divg(v+2b)|^2\divg vd\mu_g\\
&+\int_{\Omega}\nabla D_t^{r-1}D_-\divg(v+2b)\cdot\nabla D_t^{r-2}D_-\divg(v+2b) d\mu_g\\
&+\int_{\Omega}[D_t,\nabla]D_t^{r-2}D_-\divg(v+2b)\cdot\nabla D_t^{r-2}D_-\divg(v+2b)d\mu_g\\
\leq& C\Vert\nabla v\Vert_{L^\infty}E_{r}^{c+}+\int_{\Omega}\nabla D_t^{r-1}D_-\divg(v+2b)\cdot\nabla D_t^{r-2}D_-\divg(v+2b) d\mu_g.
\end{aligned}\]
Since $D_t^{r-1}D_-\divg(v+2b)$ vanishes on $\partial\Omega$, from Green's formula it follows
\[\begin{aligned}
\frac{1}{2}\frac{d}{dt}E_r^{c+}(t)\leq&C\Vert\nabla v\Vert_{L^\infty(\Omega)}E_r^{c+}-\int_{\Omega}\calt(\Delta D_t^{r-2}D_-\divg(v+2b))^2d\mu_g\\
&-\int_{\Omega}(D_t-\calt\Delta)D_t^{r-2}D_-\divg(v+2b)\cdot\Delta D_t^{r-2}D_-\divg(v+2b)d\mu_g\\
\leq&C\Vert\nabla v\Vert_{L^\infty(\Omega)}E_r^{c+}+C\epsilon_1\Vert(D_t-\calt\Delta)D_t^{r-2}D_-\divg(v+2b)\Vert_{L^2(\Omega)}^2. 
\end{aligned}\]
There is an analogous inequality for $\frac{d}{dt}E_r^{c-}$. Then (\ref{eq:dtErcest}) will be clear from following estimates: 

\begin{lemma}
(i)
\begin{equation}\label{eq:Dvbest}
\Vert D_tD_-(v+b)\Vert_{r-1,r-3}^2+\Vert D_tD_+(v-b)\Vert_{r-1,r-3}^2\leq C(\mathcal A,\mathcal E_{r-1})\mathcal E_r. 
\end{equation}
(ii)
\begin{equation}\label{eq:commTest}
\Vert D_t^{r-2}D_-[D_+,\Delta]\calt\Vert_{r-2,r-2}^2+\Vert D_t^{r-2}D_+[D_-,\Delta]\calt\Vert_{r-2,r-2}^2\leq C(\mathcal A,\mathcal E_{r-1})\mathcal E_r. 
\end{equation}
(iii) 
\begin{equation}\label{eq:Rsest}
\Vert R_{r-1}^+\Vert_{L^2(\Omega)}^2+\Vert R_{r-1}^-\Vert_{L^2(\Omega)}^2\leq C(\mathcal A,\mathcal E_{r-1})\mathcal E_r. 
\end{equation}
\end{lemma}  
\textit{Proof: }(i) Recall (\ref{eq:Dmpvb}). By (\ref{eq:intregI}), (\ref{eq:prodint}) and (\ref{eq:nonlin2}) we have
\[
\Vert D_t(b\divg(v\pm2b))\Vert_{r-1,r-3}^2\leq C(\mathcal A,\mathcal E_{r-1})\mathcal E_r. 
\] 
For $D_t(\calt\nabla q)$, commute $D_t$ with $\calt\nabla$ and specially note (\ref{eq:qreg1}). By (\ref{eq:spcomr}) and (\ref{eq:prodint}) the result thus follows. \\
(ii) Expand $D_-[\Delta,D_+]\calt$ into four terms: 
\[\begin{aligned}
D_-\Delta(v+b)^k\nabla_k\calt+\Delta(v+b)^kD_-\nabla_k\calt+2D_-\nabla(v+b)^k\cdot\nabla_k\nabla\calt\\
+2\nabla(v+b)^k\cdot D_-\nabla_k\nabla\calt. 
\end{aligned}\]
Using that $\Vert D_t^{r-2}(\phi\psi)\Vert_{L^2(\Omega)}\leq C(K)\Vert\phi\Vert_{r-2,r-2}\Vert\psi\Vert_{r-2,r-2}$ we have
\[
\Vert D_t^{r-2}(\Delta(v+b)^kD_-\nabla_k\calt)\Vert_{L^2(\Omega)}\leq C(\mathcal A,\mathcal E_{r-1})\sqrt{\mathcal E_r}\Vert D_-\nabla\calt\Vert_{r-2,r-2},
\]
\[
\Vert D_t^{r-2}(\nabla(v+b)^k\cdot D_-\nabla_k\nabla\calt)\Vert_{L^2(\Omega)}\leq C(\mathcal A,\mathcal E_{r-1})\Vert D_-\nabla^2\calt\Vert_{r-2,r-2}. 
\]
Commuting derivatives and using $D_-\calt=\calt\divg(v-2b)$, (\ref{eq:divvbr}), (\ref{eq:divvbr-1}), we can bound the right-hand side by $C(\mathcal A,\mathcal E_{r-1})\sqrt{\mathcal E_r}$. For the remaining two terms, using $\Vert D_t^{r-2}(\phi\psi)\Vert_{L^2(\Omega)}\leq C(K)(\Vert D_t\phi\Vert_{r-2,r-2}\Vert\psi\Vert_{r-3,r-3}+\Vert\phi\Vert_{r-3,r-3}\Vert D_t\psi\Vert_{r-2,r-2})$, we have
\[\begin{aligned}
\Vert D_t^{r-2}(D_-\Delta(v+b)^k\nabla_k\calt)\Vert_{L^2(\Omega)}\leq C(\mathcal A,\mathcal E_{r-1})(\Vert D_tD_-\Delta(v+b)\Vert_{r-2,r-2}\\
+\sqrt{\mathcal E_r}\Vert D_t\nabla\calt\Vert_{r-2,r-2}),
\end{aligned}\]
\[\begin{aligned}
\Vert D_t^{r-2}(D_-\nabla(v+b)^k\cdot\nabla_k\nabla\calt)\Vert_{L^2(\Omega)}\leq C(\mathcal A,\mathcal E_{r-1})(\Vert D_tD_-\nabla(v+b)\Vert_{r-2,r-2}\\
+\Vert D_t\nabla^2\calt\Vert_{r-2,r-2}). 
\end{aligned}\]
Using (\ref{eq:Dvbest}) and $D_t\calt=\calt\divg v$ we can bound the right-hand side by $C(\mathcal A,\mathcal E_{r-1})\sqrt{\mathcal E_r}$. This justifies (\ref{eq:commTest}). \\
(iii) The estimate (\ref{eq:Rsest}) simply follows from (\ref{eq:intregI}), (\ref{eq:commr1}) and (\ref{eq:prodint}). This completes the proof of the lemma as well as (\ref{eq:dtErcest}). $\Box$

\subsection{Estimates of temporal derivatives of $E_r^a,E_r^b$}\label{sec:dtErab}

In this subsection we prove
\begin{equation}\label{eq:dtEraest}
\frac{d}{dt}E_r^a(t)\leq C(\mathcal A,\mathcal E_{r-1})\mathcal E_r(t),
\end{equation}
\begin{equation}\label{eq:dtErbest}
\frac{d}{dt}E_r^b(t)\leq C(\mathcal A,\mathcal E_{r-1})\mathcal E_r(t). 
\end{equation}

\begin{lemma}
(i) It holds on $\partial\Omega$ that
\begin{equation}\label{eq:hypq}
-D_t\nabla^rq+\frac{1}{\fra}N^k\nabla^rv_k=-\nabla^rD_tq+\sum_{m=1}^{r-1}\binom{r}{m}\nabla^{r-m}v^k\widetilde\otimes\nabla_k\nabla^mq,
\end{equation}
where $\widetilde\otimes$ is a symmetrization over some of the indices.\\
(ii) Let $\curl u=\nabla\wedge u$ and $(\alpha\wedge\beta)_{ij}=\alpha_i\beta_j-\alpha_j\beta_i$. Then we have
\begin{equation}\label{eq:curlv}
\begin{aligned}
\calt^{-1}D_t\curl v-B^i\nabla_i\curl B=&\calt^{-1}[D_t,\curl]v+[\curl,B^i\nabla_i]B\\
&+\nabla\log\calt\wedge(B^i\nabla_iB)-\nabla\log\calt\wedge\nabla q,
\end{aligned}
\end{equation}
\begin{equation}\label{eq:curlB}
\begin{aligned}
D_t\curl B-B^i\nabla_i\curl v=[D_t,\curl]B+[\curl,B^i\nabla_i]v&-\nabla\divg v\wedge B\\
&-\divg v\curl B. 
\end{aligned}
\end{equation}
\end{lemma}
\textit{Proof: }Equation (\ref{eq:hypq}) is derived from the calculation
\[\begin{aligned}
\nabla^r D_tq-D_t\nabla^rq=&[\nabla^r,D_t]q=\nabla^rv^k\nabla_kq+\sum_{m=1}^{r-1}\binom{r}{m}\nabla^{r-m}v^k\nabla_k\nabla^mq
\end{aligned}\]
and the fact that $-\mathfrak a\nabla q=N$ on $\partial\Omega$. Equations (\ref{eq:curlv}) and (\ref{eq:curlB}) are derived by applying $\curl$ to the first two equations in (\ref{eq:LM}). $\Box$

\begin{lemma}
(i) Introducing the notation $Q_t(\alpha,\beta)=D_t(Q^{IJ})\alpha_I\beta_J$, $Q_{,k}(\alpha,\beta)=\nabla_k(Q^{IJ})\alpha_I\beta_J$, we have
\begin{equation}\label{eq:dtEra}
\begin{aligned}
\frac{1}{2}\frac{d}{dt}E_r^{a}(t)=&\int_{\Omega}Q(\nabla^r\divg v,\nabla^rq)+Q(\nabla^r(B^k\divg v),\nabla^rB_k)d\mu_g\\
&+\int_{\partial\Omega}Q\Big(D_t\nabla^rq-\frac{1}{\mathfrak a}N^k\nabla^rv_k,\ \nabla^rq\Big)\mathfrak ad\mu_\zeta+\mathcal I_1^a+\mathcal I_2^a+\mathcal I_3^a,
\end{aligned}
\end{equation}
where
\begin{equation}\label{eq:I1adef}
\begin{aligned}
2\mathcal I_1^a=&\int_{\Omega}\sum_{k=1}^3\Big[\calt^{-1}Q_t(\nabla^rv_k,\nabla^rv_k)dx+(Q_t+\divg v\,Q)(\nabla^rB_k,\nabla^rB_k)\Big]d\mu_g\\
&+\int_{\partial\Omega}\bigg[Q_t(\nabla^rq,\nabla^rq)+Q(\nabla^rq,\nabla^rq)\frac{D_t(\mathfrak ad\mu_\zeta)}{\mathfrak ad\mu_\zeta}\bigg]\mathfrak ad\mu_\zeta, 
\end{aligned}
\end{equation}
\begin{equation}\label{eq:I2adef}
\mathcal I_2^a=\int_{\Omega}\sum_{k=1}^n\Big[Q_{,k}(\nabla^rv_k,\nabla^rq)-B^iQ_{,i}(\nabla^rv_k,\nabla^rB_k)\Big]d\mu_g
\end{equation}
and
\begin{equation}\label{eq:I3adef}
\begin{aligned}
\mathcal I_3^a=&\int_{\Omega}\sum_{k=1}^n\Big\{Q([\calt^{-1}D_t,\nabla^r]v_k,\nabla^rv_k)+Q([D_t,\nabla^r]B_k,\nabla^rB_k)\Big\}d\mu_g\\
&+\int_{\Omega}\sum_{k=1}^n\Big\{Q([\nabla^r,B^i\nabla_i]B_k,\nabla^rv_k)+Q([\nabla^r,B^i\nabla_i]v_k,\nabla^rB_k)\Big\}d\mu_g. 
\end{aligned}
\end{equation}
(ii)
\begin{equation}\label{eq:dtErb}
\begin{aligned}
\frac{1}{2}\frac{d}{dt}E_r^{b}(t)=&\int_{\Omega}\nabla^{r-1}(\calt^{-1}D_t\curl v-B^i\nabla_i\curl B)\cdot\nabla^{r-1}\curl vd\mu_g\\
+&\int_{\Omega}\nabla^{r-1}(D_t\curl B-B^i\nabla_i\curl v)\cdot\nabla^{r-1}\curl Bd\mu_g\\
+&\frac{1}{2}\int_{\Omega}|\nabla^{r-1}\curl B|^2\divg v d\mu_g+\mathcal I^b,
\end{aligned}
\end{equation}
where
\begin{equation}\label{eq:Ibdef}
\begin{aligned}
\mathcal I^b=&\int_{\Omega}[\calt^{-1}D_t,\nabla^{r-1}]\curl v\cdot\nabla^{r-1}\curl v+[D_t,\nabla^{r-1}]\curl B\cdot\nabla^{r-1}\curl Bd\mu_g\\
&+\int_{\Omega}[\nabla^{r-1},B^i\nabla_i]\curl v\cdot\nabla^{r-1}\curl B+[\nabla^{r-1},B^i\nabla_i]\curl B\cdot\nabla^{r-1}\curl vd\mu_g.
\end{aligned}
\end{equation}
\end{lemma}
\textit{Proof: }(i) By (\ref{eq:Dtmu}) and $D_t\calt=\calt\divg v$ we have
\[\begin{aligned}
\frac{1}{2}\frac{d}{dt}\int_{\Omega}\calt^{-1}\sum_{k=1}^nQ(\nabla^rv_k,\nabla^rv_k)d\mu_g\\
=\frac{1}{2}\int_{\Omega}\calt^{-1}\sum_{k=1}^nQ_t(\nabla^rv_k,\nabla^rv_k)d\mu_g&+\int_{\Omega}\sum_{k=1}^nQ([\calt^{-1}D_t,\nabla^r]v_k,\nabla^rv_k)d\mu_g\\
&+\int_{\Omega}\sum_{k=1}^nQ(\nabla^r(\calt^{-1}D_tv_k),\nabla^rv_k)d\mu_g,
\end{aligned}\]
\[\begin{aligned}
&\frac{1}{2}\frac{d}{dt}\int_{\Omega}\sum_{k=1}^nQ(\nabla^rB_k,\nabla^rB_k)d\mu_g\\
=&\frac{1}{2}\sum_{k=1}^n\int_{\Omega}Q_t(\nabla^rB_k,\nabla^rB_k)+(\divg v)Q(\nabla^rB_k,\nabla^rB_k)d\mu_g\\
&+\int_{\Omega}\sum_{k=1}^nQ([D_t,\nabla^r]B_k,\nabla^rB_k)d\mu_g+\int_{\Omega}\sum_{k=1}^nQ(\nabla^rD_t B_k,\nabla^rB_k)d\mu_g,
\end{aligned}\]
and
\[\begin{aligned}
&\frac{1}{2}\frac{d}{dt}\int_{\partial\Omega}Q(\nabla^rq,\nabla^rq)\mathfrak ad\mu_\zeta\\
=&\frac{1}{2}\int_{\partial\Omega}\bigg[Q_t(\nabla^rq,\nabla^rq)+Q(\nabla^rq,\nabla^rq)\frac{D_t(\mathfrak ad\mu_\zeta)}{\mathfrak ad\mu_\zeta}\bigg]\mathfrak ad\mu_\zeta\\
&+\int_{\partial\Omega}Q\Big(D_t\nabla^rq-\frac{1}{\mathfrak a}N^k\nabla^rv_k,\ \nabla^rq\Big)\mathfrak ad\mu_\zeta+\int_{\partial\Omega}N^kQ(\nabla^rv_k,\nabla^rq)d\mu_\zeta.
\end{aligned}\]
By the divergence formula, we have
\[\begin{aligned}
\int_{\partial\Omega}N^kQ(\nabla^rv_k,\nabla^rq)d\mu_\zeta=\int_{\Omega}\sum_{k=1}^nQ_{,k}(\nabla^rv_k,\nabla^rq)+Q(\nabla^r\divg v,\nabla^rq)d\mu_g\\
+\int_{\Omega}\sum_{k=1}^nQ(\nabla^rv_k,\nabla^r\nabla_kq)d\mu_g. 
\end{aligned}\]
The last terms in the above identities are coupled by (\ref{eq:LM}) as
\[\begin{aligned}
\int_{\Omega}\sum_{k=1}^n\big[Q(\nabla^r(\calt^{-1}D_tv_k+\nabla_kq),\nabla^rv_k)+Q(\nabla^rD_t B_k,\nabla^rB_k)\big]d\mu_g\\
=\int_{\Omega}\sum_{k=1}^n\big\{Q([\nabla^r,B^i\nabla_i]B_k,\nabla^rv_k)+Q([\nabla^r,B^i\nabla_i]v_k,\nabla^rB_k)\big\}d\mu_g\\
+\int_{\Omega}\sum_{k=1}^n\big[Q(B^i\nabla_i\nabla^rB_k,\nabla^rv_k)+Q(\nabla^rB_k,B^i\nabla_i\nabla^rv_k)\big]d\mu_g. 
\end{aligned}\]
Again by the divergence formula, the last line equals
\[
-\int_{\Omega}B^iQ_{,i}(\nabla^rv_k,\nabla^rB_k)d\mu_g. 
\]
Combining these equations gives (\ref{eq:dtEra}). \\
(ii) Similarly, (\ref{eq:Dtmu}) and $D_t\calt=\calt\divg v$ yield 
\[\begin{aligned}
\frac{1}{2}\frac{d}{dt}\int_{\Omega}\calt^{-1}|\nabla^{r-1}\curl v|^2+|\nabla^{r-1}\curl B|^2d\mu_g-\int_{\Omega}\frac{\divg v}{2}|\nabla^{r-1}\curl B|^2d\mu_g\\
=\int_{\Omega}[\calt^{-1}D_t,\nabla^{r-1}]\curl v\cdot\nabla^{r-1}\curl B+[D_t,\nabla^{r-1}]\curl B\cdot\nabla^{r-1}\curl vd\mu_g\\
+\int_{\Omega}\nabla^{r-1}(\calt^{-1}D_t\curl v)\cdot\nabla^{r-1}\curl v+\nabla^{r-1}D_t\curl B\cdot\nabla^{r-1}\curl Bd\mu_g.
\end{aligned}\]
Writing 
\[
\calt^{-1}D_t\curl v=(\calt^{-1}D_t\curl v-B^i\nabla_i\curl B)+B^i\nabla_i\curl B,
\]
\[
D_t\curl B=(D_t\curl B-B^i\nabla_i\curl v)+B^i\nabla_i\curl v
\]
and integrating by parts, we obtain (\ref{eq:dtErb}). $\Box$

\begin{lemma}\label{lem:hypest}
\begin{equation}\label{eq:derDtq}
\Vert\Pi(\nabla^rD_tq)\Vert_{L^2(\partial\Omega)}^2+\sum_{k=1}^{r-1}\Vert \nabla^{r-k}v\nabla^{k+1}q\Vert_{L^2(\partial\Omega)}^2\leq C(\mathcal A,\mathcal E_{r-1})\mathcal E_{r}. 
\end{equation}
\end{lemma}
\textit{Proof: }From (\ref{eq:tr}) and (\ref{eq:intregI}) we see $\Vert\nabla^{r-1}D_tq\Vert_{L^2(\partial\Omega)}^2\leq C(\mathcal A,\mathcal E_{r-1})\mathcal E_{r}$. To further control $\Pi(\nabla^rD_tq)$ invoke (\ref{eq:proj}) and (\ref{eq:prodbdry}). To bound the product involving $\overline\nabla^{r_i}N,\nabla^sD_tq$, note that $\sum_{i=1}^m(r-2-r_i)+(r-1-s)\geq\frac{n-1}{2}m$ and that equality holds only when $r=\frac{n+5}{2}$, we then have
\[\begin{aligned}
\Vert\Pi(\nabla^rD_tq)\Vert_{L^2(\partial\Omega)}\leq L\Vert\overline\nabla^{r-1}N\Vert_{L^2(\partial\Omega)}+&C\left(K,\sum_{j=1}^{r-2}\Vert\overline\nabla^jN\Vert_{L^2(\partial\Omega)}\right)\\
&\times\left(\sum_{j=1}^{r-1}\Vert\nabla^jD_tq\Vert_{L^2(\partial\Omega)}+\Vert\nabla^2D_tq\Vert_{L^\infty}\right).
\end{aligned}\]
Since $\Vert\nabla^2D_tq\Vert_{L^\infty(\Omega)}\leq C(K)\Vert D_tq\Vert_r\leq C(\mathcal A,\mathcal E_{r-1})\sqrt\mathcal E_{r}$ thanks to Sobolev embedding and (\ref{eq:intregI}), we obtain the desired bound for $\Vert\Pi(\nabla^rD_tq)\Vert_{L^2(\partial\Omega)}$.  \\
For the other term, we see from (\ref{eq:tr}), Sobolev embedding and (\ref{eq:intregI}) that $\Vert\nabla^jv\Vert_{L^2(\partial\Omega)}+\Vert\nabla^2v\Vert_{L^\infty}\leq C(\mathcal A,\mathcal E_{r-1})\sqrt{\mathcal E_r}$, $j\leq r-1$. Using $\Vert\nabla^jq\Vert_{L^2(\partial\Omega)}\leq C(\mathcal A,\mathcal E_{r-1})$, $j\leq r-1$ and $|\nabla^2q|\leq L$, we obtain the bound for $\nabla^{r-k}v\nabla^{k+1}q$, $k\leq r-2$. Finally, for $k=r-1$ we use (\ref{eq:qNbdry}) and $\Vert\nabla v\Vert_{L^\infty}\leq C(K)\sqrt{\mathcal E_{r_*}}$. $\Box$\\[2mm]
\textit{Proof of (\ref{eq:dtEraest}): }We make use of (\ref{eq:dtEra}). By (\ref{eq:intregI}) and (\ref{eq:prodint}), we have
\[
\left|\int_{\Omega}Q(\nabla^r\divg v,\nabla^rq)+Q(\nabla^r(B^k\divg v),\nabla^rB_k)d\mu_g\right|\leq C(\mathcal A,\mathcal E_{r-1})\mathcal E_{r}. 
\]
Moreover, by (\ref{eq:derDtq}) and (\ref{eq:Era}), we have
\[\begin{aligned}
\int_{\partial\Omega}Q\Big(D_t\nabla^rq-\frac{1}{\mathfrak a}N^k\nabla^rv_k,\ \nabla^rq\Big)\mathfrak ad\mu_\zeta\leq& \Vert D_t\nabla^rq-\frac{1}{\mathfrak a}N^k\nabla^rv_k\Vert_{L^2(\partial\Omega)}\Vert\Pi(\nabla^rq)\Vert_{L^2(\partial\Omega)}\\
\leq&C(\mathcal A,\mathcal E_{r-1})\mathcal E_{r}. 
\end{aligned}\]
For $\mathcal I_1^a$ and $\mathcal I_2^a$, it follows from (\ref{eq:Dtmu})-(\ref{eq:DtNbdry}), (\ref{eq:ass2}) and the calculation
\begin{equation}\label{eq:Dtloga}
|D_t\log\mathfrak a|=|\nabla q|^{-1}\big|D_t|\nabla q|\big|\leq |\nabla q|^{-1}|D_t\nabla q|=|\nabla q|^{-1}|\nabla v^k\nabla_kq+\nabla D_tq|
\end{equation}
that $|\mathcal I_1^a|+|\mathcal I_2^a|\leq C(\mathcal A,\mathcal E_{r-1})\mathcal E_{r}$. Finally, from (\ref{eq:commr1}) it follows that $|\mathcal I_3^a|\leq C(\mathcal A,\mathcal E_{r-1})\mathcal E_{r}$. This completes the proof. $\Box$\\[2mm]
\textit{Proof of (\ref{eq:dtErbest}): }The result simply follows from (\ref{eq:curlv})-(\ref{eq:curlB}), (\ref{eq:dtErb})-(\ref{eq:Ibdef}) and (\ref{eq:intregI}), (\ref{eq:commr1}), (\ref{eq:prodint}). $\Box$

\subsection{Completion of the proof}\label{sec:pfmainfin}

In the previous two subsections we have shown that $\frac{d}{dt}E_r(t)\leq C(\mathcal A,\mathcal E_{r-1})\mathcal E_r(t)$ for $r>r_*$. A much easier calculation shows that $\frac{d}{dt}\mathcal E_{r_*}(t)\leq C(\mathcal A,\mathcal E_{r_*})\mathcal E_{r_*+1}(t)$. Thus (\ref{eq:dtErest}) is proven, under the assumptions (\ref{eq:infT})-(\ref{eq:ass2}), which we shall recover now. Setting $(\epsilon_0,\epsilon_1,K,L)=(\frac{1}{2}\epsilon_{0*},\epsilon_{1*},2K_*,2L_*)$, we claim that (\ref{eq:recTS})-(\ref{eq:volest}) hold true on $[0,T]$ for some definite $T>0$. Indeed, (\ref{eq:recinfT}) follows the maximal principle for the parabolic equation $D_t\calt=\calt\Delta\calt$. For (\ref{eq:recTS}), (\ref{eq:recass2}) and (\ref{eq:volest}) we control the material derivative of the variables on the left-hand side: 
\[
|D_t\nabla q|+|D_t\nabla D_tq|\leq C(\mathcal A,\mathcal E_{r_*})\sqrt{\mathcal E_{r_*+1}},
\]
\[
|D_t\nabla^2q|\leq|[D_t,\nabla^2]q|+|\nabla^2D_tq|\leq C(\mathcal A,\mathcal E_{r_*})\sqrt{\mathcal E_{r_*+1}},
\]
\[
\left|\frac{d}{dt}\log(\vol\Omega)\right|\leq\Vert\divg v\Vert_{L^\infty(\Omega)}\leq C(\mathcal A,\mathcal E_{r_*}),
\]
thanks to (\ref{eq:intregI}) and Sobolev embedding. For (\ref{eq:recass1}), we invoke (3.86) from \cite{LZeng}:
\begin{equation}\label{eq:thetaptbd}
|D_t\theta|\leq |\nabla^2v|+C|\theta||\nabla v|,
\end{equation}
which can be controlled by $C(\mathcal A,\mathcal E_{r_*})\sqrt{\mathcal E_{r_*+1}}$. It remains to estimate the injectivity radius at time $t$. For this, it is more convenient to work with Eulerian coordinates. We need the following lemmas, the proofs of which are delayed to the end of this section. 

\begin{lemma}\label{lem:thetaest}
Let $\Sigma\subset\mathbb R^n$ be a closed hypersurface with bounded second fundamental form $\theta$ and positive injectivity radius $\iota_0$ of normal exponential map. Then we have
\[
|\theta|\leq\frac{1}{\iota_0}. 
\]
\end{lemma}

\begin{lemma}\label{lem:injradest}
Let $\Sigma,\Sigma_t\subset\mathbb R^n$ be closed hypersurfaces with bounded second fundamental forms $\theta,\theta_t$. Suppose that
\begin{equation}\label{eq:dLest}
\begin{aligned}
d_L(\Sigma,\Sigma_t):=\sup\{|\bar x_t-\bar x|+|N_t(\bar x_t)-N(\bar x)|:\bar x_t\in\Sigma_t,\bar x\in\Sigma,|\bar x_t-\bar x|=d(\bar x_t,\Sigma)\}\\
\leq Mt
\end{aligned}
\end{equation}
and that
\begin{equation}\label{eq:thetaerr}
|\theta_t|\leq\frac{1}{\iota_0(\Sigma)}+\Lambda t. 
\end{equation}
Then for every $\epsilon\in(0,\iota_0(\Sigma))$, there is a $t_0=t_0(\epsilon,\iota_0(\Sigma),M,\Lambda)>0$ such that $\iota_0(\Sigma_t)\geq\iota_0(\Sigma)-\epsilon$, $\forall t\in(0,t_0]$. 
\end{lemma}
With $\Sigma=\partial\mathcal D_0$ and $\Sigma_t=\partial\mathcal D_t$, condition (\ref{eq:thetaerr}) is satisfied thanks to Lemma \ref{lem:thetaest} and (\ref{eq:thetaptbd}). To verify (\ref{eq:dLest}), one considers the flow map $x(t,y):y\in\partial\mathcal D_0\mapsto x\in\partial\mathcal D_t$. Since $\Vert \partial_t|_{y=\text{const}}x(t,y)\Vert_{C^1_y}=\Vert v(t,x(t,y))\Vert_{C^1_y}$ and $\Vert\theta\Vert_{L^\infty(\partial\mathcal D_0)}$ are uniformly bounded, (\ref{eq:dLest}) thus follows. Till now the proof of Theorem \ref{thm:main} is completed. \\[2mm]
\textit{Proof of Lemma \ref{lem:thetaest}: }Suppose $|\theta(\bar x_0)|>\frac{1}{\iota_0}$, then for some $s_0\in(-\iota_0,\iota_0)$ the matrix $I_{n-1}+s_0\theta(\bar x_0)$ has negative eigenvalues. Then consider the function $h(\bar x)=\frac{1}{2}|\bar x-x_0|^2$ where $x_0=\bar x_0+s_0N(\bar x_0)$. A calculation using the graph coordinates introduced in Lemma \ref{lem:grcoor} gives $\partial^2h(\bar x_0)=I_{n-1}+s_0\theta(\bar x_0)$, hence $h$ is not minimizing at $\bar x_0$, which contradicts the definition of $\iota_0$. $\Box$\\[2mm]
\textit{Proof of Lemma \ref{lem:injradest}: }Let $\bar x_t\in\Sigma_t$, $\bar x\in\Sigma$ be such that $|\bar x-\bar x_t|=d(\bar x_t,\Sigma)$. Let $x=\bar x+\iota_0(\Sigma)N(\bar x)$, $r=2d_L(\Sigma,\Sigma_t)$ and $R=\iota_0(\Sigma)-r$, then $B_R(x)\cap\Sigma_t=\emptyset$. With $R_t=R-t^\frac{1}{5}$, $x_t=\bar x_t+R_tN_t(\bar x_t)$, we shall show that
\begin{equation}\label{eq:Btnointsec}
B_{R_t}(x_t)\cap\Sigma_t=\{\bar x_t\}
\end{equation}
for $t\ll1$. As a first step, we estimate the set $B_{R_t}(x_t)\backslash B_R(x)$. Suppose $z$ lies in the set, that is, 
\[
|z-x|^2>R^2,\qquad |z-x_t|^2\leq R_t^2,
\]
then we immediately get
\[
R<|z-x|\leq R_t+|x-x_t|,
\]
\[
\frac{z-x}{|z-x|}\cdot\frac{x_t-x}{|x_t-x|}>\frac{R-R_t}{|x_t-x|}\frac{R+R_t}{2|z-x|}+\frac{|x-x_t|}{2|z-x|}.
\]
Noting $R-R_t=t^\frac{1}{5}$, $r\leq 2Mt$ and
\[
|x_t-x|=|\bar x-\bar x_t+(R+r)N(\bar x)-R_tN_t(\bar x_t)|=t^\frac{1}{5}+O(r),
\]
we are leading to
\[
\cos\phi=\frac{z-x}{|z-x|}\cdot\frac{x_t-x}{|x_t-x|}>1-Ct^\frac{4}{5},
\]
where $\phi$ is the angle formed by $z-x$ and $x_t-x$. Combined with the estimate of $|z-x|$, we find that 
\[
B_{R_t}(x_t)\backslash B_R(x)\subset B_{Ct^\frac{2}{5}}(\bar x_t), 
\]
where $C$ only depends on $\iota_0(\Sigma)$ and $M$. For clarity let us drop the constant $C$ and simply write $B_{t^\frac{2}{5}}(\bar x_t)$. 

Now that (\ref{eq:Btnointsec}) is reduced to the local statement
\begin{equation}\label{eq:locnointsec}
B_{R_t}(x_t)\cap\left(B_{t^\frac{2}{5}}(\bar x_t)\cap\Sigma_t\right)=\{\bar x_t\},
\end{equation}
we invoke Lemma \ref{lem:grcoor}. It states that the graph parametrization $\mathbf f(\xi)=(\xi,f(\xi))$ around $\bar x_t$ is defined for $|\xi|\leq 2t^\frac{2}{5}$ and satisfies the estimates (\ref{eq:graphest}) with $\epsilon\leq CKt^\frac{2}{5}$ and $K=\frac{2}{\iota_0(\Sigma)}>|\theta_t|$. We claim that
\[
B_{t^\frac{2}{5}}(\bar x_t)\cap\Sigma_t\ \subset\ S:=\{\mathbf f(\xi)\mid|\xi|\leq 2t^\frac{2}{5}\}, 
\]
or in other words, the projection $\pi_t:\mathbb R^n\to\mathbf T_{\bar x_t}\Sigma_t$ is injective on $B_{t^\frac{2}{5}}(\bar x_t)\cap\Sigma_t$. 
Since $\theta_t$ is uniformly bounded, we only need to exclude that there are multiple sheets of $\Sigma_t$ in $B_{t^\frac{2}{5}}(\bar x_t)$. This is clear from $\iota_0(\Sigma)>0$ and (\ref{eq:dLest}), and the claim thus follows. To prove (\ref{eq:locnointsec}), it suffices to show that $\xi=0$ uniquely minimizes $h(\xi)=\frac{1}{2}|\mathbf f(\xi)-x_t|^2$ in $\{|\xi|\leq2t^\frac{2}{5}\}$. A calculation gives 
\[
\partial_\alpha\partial_\beta h=\delta_{\alpha\beta}+\partial_\beta((f-R_t)\partial_\alpha f). 
\]
Now that $|\partial f|=O(t^\frac{2}{5})$, $|R_t-f|=\iota_0(\Sigma)-t^\frac{1}{5}+O(t^\frac{2}{5})$ and $\sum|\partial_\alpha\partial_\beta f|^2\leq(\frac{1}{\iota_0(\Sigma)}+\Lambda t)^2+O(t^\frac{2}{5})$, we have $\sum|\partial_\beta((f-R_t)\partial_\alpha f)|^2<1$, hence $\partial^2h$ is positive definite for small $t$. Suppose $\xi\neq0$ minimizes $h$, then the function $H(s)=h(s\xi)$ satisfies $H'(0)=0=H'(1)$ and $H''(s)>0$, $s\in[0,1]$, which is impossible. This justifies (\ref{eq:locnointsec}) and (\ref{eq:Btnointsec}). A similar argument shows that
\[
B_{R_t}(x'_t)\cap\Sigma_t=\{\bar x_t\},\qquad x'_t=\bar x_t-R_tN_t(\bar x_t), 
\]
and thus $\iota_0(\Sigma_t)\geq R_t$. This completes the proof. $\Box$

\begin{lemma}\label{lem:grcoor}
Let $\Sigma\subset\mathbb R^n$ be a closed hypersurface with $|\theta|\leq K$. Suppose $0\in\bar x_0\in\Sigma$, $N(\bar x_0)=(0,...,0,1)$ and write $x=(\xi,x^n)$, $\xi\in\mathbb R^{n-1}=\mathbf T_{\bar x_0}\Sigma$. Then for every $\epsilon>0$, there is a ball $B^{n-1}_{r_0}(0)\subset\mathbf T_{\bar x_0}\Sigma$ of radius $r_0>\frac{\epsilon}{CK}$ and a function $f:B^{n-1}_{r_0}\to\mathbb R$ such that $\xi\mapsto\mathbf f(\xi)=(\xi,f(\xi))$ is a diffeomorphism onto an open neighborhood in $\Sigma$ and that
\begin{equation}\label{eq:graphest}
f(0)=0,\quad\partial_\alpha f(0)=0,\quad\sum|\partial_\alpha\partial_\beta f|^2\leq(1+\epsilon)|\theta|^2, 
\end{equation}
where $\partial_\alpha=\frac{\partial}{\partial\xi^\alpha}$ and $\alpha,\beta$ range from $1$ to $n-1$. As a consequence, we have
\[
\mathbf f(B^{n-1}_{r_0})\supset B_{r_0}^\Sigma(\bar x_0),
\]
\[
d_\Sigma(\mathbf f(\xi),\mathbf f(\eta))\leq2|\xi-\eta|,\quad\forall\xi,\eta\in B^{n-1}_{r_0},
\]
where $d_\Sigma$ is the geodesic distance and $B_{r_0}^\Sigma$ is the geodesic ball. 
\end{lemma}
\textit{Proof: }By implicit function theorem $f$ locally exists. A calculation gives
\[
N_\alpha=\frac{-\partial_\alpha f}{\sqrt{1+|\partial f|^2}},\quad N_n=\frac{1}{\sqrt{1+|\partial f|^2}},\quad\zeta_{\alpha\beta}=\delta_{\alpha\beta}-\partial_\alpha f\partial_\beta f
\]
where $|\partial f|^2=\delta^{\alpha\beta}\partial_\alpha f\partial_\beta f$ and $\zeta_{\alpha\beta}$ is the induced metric on $\Sigma$. Then we have
\[
\zeta^{\alpha\beta}=\delta^{\alpha\beta}+O(|\partial f|^2),
\]
\[
|\theta|^2=\zeta^{\alpha\beta}\partial_\alpha N_i\partial_\beta N_j\delta^{ij}=\frac{|\partial^2f|^2}{\sqrt{1+|\partial f|^2}}+O(|\partial f||\partial^2f|^2),
\]
\[
\text{Length of }t\mapsto\mathbf f(t\xi)=\int_0^1(|\xi|^2-\xi^\alpha\xi^\beta\partial_\alpha f\partial_\beta f)^\frac{1}{2}dt\leq|\xi|^2(1+\frac{1}{2}\Vert\partial f\Vert_{L^\infty}^2).
\]
Now the lemma follows a bootstrap argument. $\Box$

\section{Regularity estimates I}\label{sec:regI}

\begin{proposition}\label{prop:regI}
Under the assumptions (\ref{eq:TS})-(\ref{eq:ass2}) and with the notation (\ref{eq:calAdef}), for $r\geq r_*$ we have
\begin{equation}\label{eq:qNbdry}
\sum_{k=1}^{r-1}\Vert\overline{\nabla}^kN\Vert_{L^2(\partial\Omega)}^2+\sum_{k=1}^{r}\Vert\overline{\nabla}^kq\Vert_{L^2(\partial\Omega)}^2\leq C(\mathcal A,\mathcal E_{r-1})\mathcal E_r(t),
\end{equation}
\begin{equation}\label{eq:intregI}
\Vert(v,B,\calt,\divg v)\Vert_{r,r-1}^2+\Vert q\Vert_{r,r-2}^2+\Vert D_tq\Vert_{r,r-3}^2\leq C(\mathcal A,\mathcal E_{r-1})\mathcal E_r(t). 
\end{equation}
\end{proposition}
The proof involves an induction on $r\geq r_*$. For $r=r_*$, by (\ref{eq:Er*def}) it remains to show $\Vert\overline{\nabla}^kN\Vert_{L^2(\partial\Omega)}^2\leq C(\mathcal A,\mathcal E_{r_*})$, $k\leq r_*-1$. For this, we invoke (\ref{eq:DNexpr}) with $s=r_*-1$, $f=q$ and (\ref{eq:prodbdry}) with $r_i=r_*-1-s_i$, $f_i=\nabla^{s_i+1}q$. Since $\sum(r_*-1-s_i)=(m-1)(r_*-1)>\frac{n-1}{2}(m-1)$, the result thus follows. 

Next, we prove (\ref{eq:qNbdry}), (\ref{eq:intregI}) for $r>r_*$ under the inductive hypothesis
\begin{equation}\label{eq:qNhypo}
\sum_{k=1}^{r-2}\Vert\overline{\nabla}^kN\Vert_{L^2(\partial\Omega)}^2+\sum_{k=1}^{r-1}\Vert\overline{\nabla}^kq\Vert_{L^2(\partial\Omega)}^2\leq C(\mathcal A,\mathcal E_{r-1}),
\end{equation}
\begin{equation}\label{eq:inthypo}
\Vert(v,B,\calt,\divg v)\Vert_{r-1,r-2}^2+\Vert q\Vert_{r-1,r-3}^2+\Vert D_tq\Vert_{r-1,r-4}^2\leq C(\mathcal A,\mathcal E_{r-1}). 
\end{equation}
The proofs are presented in the following three subsections. 

\subsection{Boundary regularity}\label{sub:bdryreg}

In this subsection we prove
\begin{equation}\label{eq:Nrec}
\sum_{k=1}^{r-1}\Vert\overline{\nabla}^kN\Vert_{L^2(\partial\Omega)}^2\leq C(\mathcal A,\mathcal E_{r-1})\mathcal E_r(t). 
\end{equation}
The proof relies on (\ref{eq:proj}) with $f=q$. By (\ref{eq:Era}), it remains to control the product $\overline\nabla^{r_1}N\cdots\overline\nabla^{r_m}N\nabla^sq$. For this, invoke (\ref{eq:prodbdry}) with $m$ replaced by $m+1$ and $f_i=\overline\nabla^{r_i}N$, $f_{m+1}=\nabla^sq$. Observe that $\sum_{i=1}^m(r-2-r_i)+r-1-s\geq\frac{n-1}{2}m$ and equality holds only when $r=\frac{n+5}{2}$, we thus have
\[\begin{aligned}
\Vert\prod_{i=1}^m\overline\nabla^{r_i}N\nabla^sq\Vert_{L^2(\partial\Omega)}\leq C(K)\left(\sum_{k=0}^{r-2}\Vert\overline\nabla^kN\Vert_{L^2(\partial\Omega)}+\Vert\overline\nabla N\Vert_{L^\infty}\right)\\
\times\left(\sum_{k=0}^{r-1}\Vert\nabla^kq\Vert_{L^2(\partial\Omega)}+\Vert\nabla^2q\Vert_{L^\infty}\right). 
\end{aligned}\]
Then (\ref{eq:Nrec}) follows from (\ref{eq:qNhypo}) and (\ref{eq:ass1}), (\ref{eq:ass2}). 

As an application, we prove the following elliptic estimate: 

\begin{proposition}\label{prop:diri2}
For $u\in H^r(\Omega)$, $u|_{\partial\Omega}=0$ we have
\begin{equation}\label{eq:diri2}
\Vert\nabla^ru\Vert_{L^2(\Omega)}^2\leq C\Vert\nabla^{r-2}\Delta u\Vert^2_{L^2(\Omega)}+C(\mathcal A,\mathcal E_{r-1})(\Vert u\Vert_{r-1}^2+\Vert\nabla u\Vert_{L^\infty(\partial\Omega)}^2\mathcal E_r)
\end{equation}
Moreover, for $2\leq s\leq r-1$ we have
\begin{equation}\label{eq:diri3}
\Vert\nabla^su\Vert_{L^2(\Omega)}^2\leq C\Vert\nabla^{s-2}\Delta u\Vert^2_{L^2(\Omega)}+C(\mathcal A,\mathcal E_{r-1})\Vert u\Vert_{s-1}^2.
\end{equation} 
\end{proposition}
\textit{Proof: }The key ingredients are (\ref{eq:diri}) and (\ref{eq:proj}). For $\Pi(\nabla^ru)$, the above arguments combined with (\ref{eq:qNhypo}), (\ref{eq:Nrec}) imply that
\[\begin{aligned}
\Vert\Pi(\nabla^ru)\Vert_{L^2(\partial\Omega)}\leq C(\mathcal A,\mathcal E_{r-1})\bigg(\Vert\nabla_Nu\Vert_{L^\infty(\partial\Omega)}\sqrt{\mathcal E_r}+\Vert\nabla^2u\Vert_{L^\infty(\partial\Omega)}\\
+\sum_{k=0}^{r-1}\Vert\nabla^ku\Vert_{L^2(\partial\Omega)}\bigg).
\end{aligned}\]
By (\ref{eq:diri}), (\ref{eq:tr}), $\Vert\nabla^2u\Vert_{L^\infty(\Omega)}\leq C(K)\Vert u\Vert_r$ and Young's inequality, (\ref{eq:diri2}) thus follows. 

To prove (\ref{eq:diri3}), rewrite (\ref{eq:proj}) as
\[
\Pi(\nabla^su)=\sum_{s_1+\cdots+s_{m}=s}c_{s_1...s_{m}}(g)\overline{\nabla}^{s_1}N\cdots\overline{\nabla}^{s_{m-1}}N\nabla^{s_m} u,
\]
where the summation is taken over $m\geq2$ and $s_i\geq1$. To control the product term, invoke (\ref{eq:prodbdry}). Observe that $\sum_{i=1}^{m-1}(r-2-s_i)+(s-1-s_m)\geq\frac{n-1}{2}(m-1)$ and equality holds only when $r=\frac{n+5}{2}$, we thus have
\[
\Vert\Pi(\nabla^ru)\Vert_{L^2(\partial\Omega)}\leq C(\mathcal A,\mathcal E_{r-1})\left(\sum_{k=0}^{s-1}\Vert\nabla^ku\Vert_{L^2(\partial\Omega)}+\sigma\Vert\nabla u\Vert_{L^\infty(\partial\Omega)}\right)
\]
where $\sigma=1$ if $s=r-1=\frac{n+3}{2}$ and $\sigma=0$ otherwise. Then (\ref{eq:diri3}) follows from (\ref{eq:diri}), (\ref{eq:tr}), $\Vert\nabla u\Vert_{L^\infty(\Omega)}\leq C(K)\Vert u\Vert_{\frac{n+3}{2}}$ and Young's inequality. $\Box$

\subsection{Estimates of $(v,B,q,\calt)$}

As a first step, we show that 
\begin{equation}\label{eq:Trr-1}
\Vert \calt\Vert_{r,r-1}^2\leq C(\mathcal A,\mathcal E_{r-1})\mathcal E_r, 
\end{equation}
\begin{equation}\label{eq:Dtr-1divvL2}
\Vert D_t^{r-1}\divg v\Vert_{L^2(\Omega)}^2\leq C(\mathcal A,\mathcal E_{r-1})\mathcal E_r. 
\end{equation}
For (\ref{eq:Trr-1}), we use $\Delta\calt=\divg v$, (\ref{eq:diri2}), $D_t\calt=\calt\divg v$, (\ref{eq:rsequiv}) and (\ref{eq:inthypo}). For (\ref{eq:Dtr-1divvL2}), recall $b=\calt^\frac{1}{2}B$ and note the formulas
\begin{equation}\label{eq:divb}
\divg b=\frac{1}{2}\calt^{-1}b^k\nabla_k\calt,
\end{equation}
\begin{equation}\label{eq:Dtdivvrec}
D_-\divg(v+2b)+D_+\divg(v-2b)=2D_t\divg v-4b^k\nabla_k\divg b,
\end{equation}
\begin{equation}\label{eq:Dtdivbrec}
D_-\divg(v+2b)-D_+\divg(v-2b)=4D_t\divg b-2b^k\nabla_k\divg v.
\end{equation}
From (\ref{eq:divb}), (\ref{eq:Trr-1}), (\ref{eq:inthypo}) and (\ref{eq:prodint}), (\ref{eq:nonlin2}) it follows that
\[
\Vert D_t^{r-2}b^k\nabla_k\divg b\Vert_{L^2(\Omega)}^2\leq C(\mathcal A,\mathcal E_{r-1})\mathcal E_r,
\]
then (\ref{eq:Dtr-1divvL2}) is clear from (\ref{eq:Dtdivvrec}), (\ref{eq:Erc}) and Poincar\'e inequality. 

\begin{lemma}\label{lem:homhd}
For $s\geq0$ we have
\begin{equation}\label{eq:dtvb}
\left\{
\begin{aligned}
&D_t^{s+1}v=D_t^s(\calt B^k\nabla_k B)-D_t^s(\calt\nabla q),\\
&D_t^{s+1}B=D_t^s(B^k\nabla_k v)-D_t^s(B\divg v),
\end{aligned}
\right.
\end{equation}
%\begin{equation}\label{eq:ellipdtsq}
%-\Delta D_t^sq=[D_t^s,\Delta]q-D_t^s(\nabla_i B^j\nabla_j B^i)+D_t^s[\divg,\calt^{-1}D_t]v+D_t^s(\calt^{-1}D_t\divg v), 
%\end{equation}
%and
\begin{equation}\label{eq:ellipdtsq2}
-\divg(\calt\nabla D_t^sq)=D_t^{s+1}\divg v+[\nabla_i,D_t^{s+1}]v^i+\divg[D_t^s,\calt\nabla]q-[\nabla_i,D_t^s(\calt B^j\nabla_j)]B^i
\end{equation}
\end{lemma}
\textit{Proof: }Equation (\ref{eq:dtvb}) is obtained by applying $D_t^s$ to the first two equations in (\ref{eq:LM}). Equation (\ref{eq:ellipdtsq2}) is obtained by taking divergence on both sides of (\ref{eq:dtvb}). $\Box$\\[2mm]
\begin{lemma}
\begin{equation}\label{eq:qreg}
\Vert q\Vert_{r,r-2}^2\leq C(\mathcal A,\mathcal E_{r-1})\mathcal E_r,
\end{equation}
\begin{equation}\label{eq:vBreg}
\Vert(v,B)\Vert_{r,r-1}^2\leq C(\mathcal A,\mathcal E_{r-1})\mathcal E_r. 
\end{equation}
\end{lemma}
\textit{Proof: }By (\ref{eq:commr}), (\ref{eq:spcomr}) and (\ref{eq:inthypo}), (\ref{eq:Dtr-1divvL2}) we have
\[
\Vert\divg(\calt\nabla D_t^sq)\Vert_{r-s-2}^2\leq C(\mathcal A,\mathcal E_{r-1})\mathcal E_r,\quad s\leq r-2. 
\]
For $s=r-2$ we first use (\ref{eq:W1pest}) and then use (\ref{eq:W2pest}), and for $s<r-2$ we directly use (\ref{eq:diri2}), (\ref{eq:diri3}), leading to 
\[
\Vert D_t^sq\Vert_{r-s}^2\leq C(\mathcal A,\mathcal E_{r-1})\mathcal E_r,\quad s\leq r-2. 
\]
This combined with (\ref{eq:rsequiv}) justifies (\ref{eq:qreg}). Next, using
\begin{equation}\label{eq:ptHodge}
|\nabla^ku|^2\leq C(|\nabla^{k-1}\divg u|^2+|\nabla^{k-1}\curl u|^2+\delta^{ij}Q(\nabla^ku_i,\nabla^ku_j))
\end{equation}
and (\ref{eq:inthypo}) we obtain $\Vert(v,B)\Vert_{r}^2\leq C(\mathcal A,\mathcal E_{r-1})\mathcal E_r$. Then (\ref{eq:vBreg}) follows (\ref{eq:dtvb}), (\ref{eq:rsequiv}) and an induction on the number of $D_t$'s. $\Box$

\subsection{Estimates of $\divg v$ and further regularity}

\begin{lemma}
\begin{equation}\label{eq:divvbr}
\Vert\divg v\Vert_{r,r-1}^2+\Vert\divg b\Vert_{r,r-1}^2\leq C(\mathcal A,\mathcal E_{r-1})\mathcal E_r. 
\end{equation}
\end{lemma}
\textit{Proof: }From (\ref{eq:Dtdivvb}), (\ref{eq:commr}), (\ref{eq:prodint}) and (\ref{eq:inthypo}) it follows that $\Vert\Delta\divg(v\pm2b)\Vert_{r-3}\leq C(\mathcal A,\mathcal E_{r-1})$. Then by (\ref{eq:diri3}) we have $\Vert\divg(v\pm2b)\Vert_{r-1}\leq C(\mathcal A,\mathcal E_{r-1})$ and by Sobolev embedding we have $\Vert\nabla\divg(v\pm2b)\Vert_{L^\infty(\Omega)}\leq C(\mathcal A,\mathcal E_{r-1})$. Further using (\ref{eq:commr1}) and (\ref{eq:vBreg}) we derive $\Vert\Delta\divg(v\pm2b)\Vert_{r-2}^2\leq C(\mathcal A,\mathcal E_{r-1})\mathcal E_r$, hence by (\ref{eq:diri2}) we have $\Vert\divg(v\pm2b)\Vert_{r}^2\leq C(\mathcal A,\mathcal E_{r-1})\mathcal E_r$. 

To recover the full regularity, we use (\ref{eq:Dtdivvb1}). By similar arguments we have
\[
\Vert D_t^{s-1}D_-\divg(v+2b)\Vert_{r-s}^2+\Vert D_t^{s-1}D_+\divg(v-2b)\Vert_{r-s}^2\leq C(\mathcal A,\mathcal E_{r-1})\mathcal E_r,\quad s=1,...,r-2. 
\]
When $s=r-1$, the inequality simply follows (\ref{eq:Erc}). By (\ref{eq:Dtdivvrec}), (\ref{eq:Dtdivbrec}) and an induction on $s=1,...,r-1$, we derive 
\[
\Vert D_t^s\divg v\Vert_{r-s}^2+\Vert D_t^s\divg b\Vert_{r-s}^2\leq C(\mathcal A,\mathcal E_{r-1})\mathcal E_r,\quad s=1,...,r-1. 
\]
Now (\ref{eq:divvbr}) follows from (\ref{eq:rsequiv}). $\Box$\\[2mm]
\textbf{Remark. }We have shown $\Vert\divg(v\pm2b)\Vert_{r-1}\leq C(\mathcal A,\mathcal E_{r-1})$ above. In fact, we can further show that
\begin{equation}\label{eq:divvbr-1}
\Vert\divg(v+2b)\Vert_{r-1,r-2}+\Vert\divg(v-2b)\Vert_{r-1,r-2}\leq C(\mathcal A,\mathcal E_{r-1}). 
\end{equation}
Indeed, by definition we have $\Vert\divg(v\pm2b)\Vert_{r-1,r-2}\leq\Vert\divg(v\pm2b)\Vert_{r-1}+\Vert D_t\divg(v\pm2b)\Vert_{r-2,r-3}$. Noting $\divg(v\pm2b)=\calt^{-1}D_\pm\calt$, commuting $D_t$ with $\calt^{-1}D_\pm$ and using $D_t\calt=\calt\divg v$, we can bound $\Vert D_t\divg(v\pm2b)\Vert_{r-2,r-3}$ by $C(\mathcal A,\mathcal E_{r-1})$ as well. 

\begin{lemma}
\begin{equation}\label{eq:qreg1}
\Vert\nabla^rq\Vert_{L^2(\partial\Omega)}^2+\Vert D_tq\Vert_{r,r-3}^2\leq C(\mathcal A,\mathcal E_{r-1})\mathcal E_r. 
\end{equation}
\end{lemma}
\textit{Proof: }From (\ref{eq:ellipdtsq2}) and (\ref{eq:commr1}), (\ref{eq:qreg}), (\ref{eq:vBreg}), (\ref{eq:divvbr}) we see 
\begin{equation}\label{eq:divTdDtsqest}
\Vert\divg(\calt\nabla D_t^sq)\Vert_{r-1-s}^2\leq C(\mathcal A,\mathcal E_{r-1})\mathcal E_r,\quad s=0,...,r-2. 
\end{equation}
Similar bounds hold for $\Delta D_t^sq$, thanks to $\divg(\calt\nabla f)=\nabla\calt\cdot\nabla f+\calt\Delta f$ and (\ref{eq:infT}). Then the bound for $\Vert D_tq\Vert_{r,r-3}$ follows from (\ref{eq:diri2}), (\ref{eq:diri3}), (\ref{eq:ass2}) and (\ref{eq:rsequiv}), and the bound for $\Vert\nabla^rq\Vert_{L^2(\partial\Omega)}$ follows from (\ref{eq:bdryHr}). $\Box$\\[2mm]
Till now the proof of (\ref{eq:qNbdry}) and (\ref{eq:intregI}) is completed.

\section{Proof of Theorem \ref{thm:buc}}

For clarity, let us drop the $t$-dependence in (\ref{eq:Mdef}), (\ref{eq:Kdef}) and simply assume
\[
\Vert\divg v\Vert_{L^\infty(\Omega)}+\Vert\curl v\Vert_{L^\infty(\Omega)}+\Vert\curl B\Vert_{L^\infty(\Omega)}+\Vert\overline\nabla(v,B)\Vert_{L^\infty(\partial\Omega)}\leq M,
\]
\[
\Vert\theta\Vert_{L^\infty(\partial\Omega)}+\frac{1}{\iota_0(\Omega)}\leq K
\]
for some constants $M,K$ and for all $t\in[0,T)$. Also, let us group all the control parameters together as 
\[
\mathcal K=(K,M,\epsilon_0,L_0,L_1). 
\]
We need lower order energies in analogy to (\ref{eq:Era})-(\ref{eq:Erc}). The definition of $E_1,E_2$ and $E_3^a,E_3^b$ have exactly the same form. However, for $E_3^{c\pm}$ we make slight modifications: 
\begin{equation}\label{eq:E3c+}
E_3^{c+}(t)=\int_\Omega |\nabla D_+D_-\divg(v+2b)|^2d\mu_g,
\end{equation}
\begin{equation}\label{eq:E3c-}
E_3^{c-}(t)=\int_\Omega |\nabla D_-D_+\divg(v-2b)|^2d\mu_g. 
\end{equation}
Further denote
\begin{equation}\label{eq:calEdef}
\mathcal E=E_1+E_2+E_3^a+E_3^b+E_3^{c+}+E_3^{c-}+2. 
\end{equation}
The constant $2$ is added so that $\log\mathcal E\geq\log2>0$.

\subsection{An integral inequality for $E_3^{c\pm}$}

We shall show that
\begin{equation}\label{eq:intE3c}
\sup_{s\leq t}E_3^{c+}(s)\leq \frac{1}{4}\sup_{s\leq t}\mathcal E(s)+C(\mathcal K)\left(\mathcal E(0)+\int_0^t\mathcal E(s)\log\mathcal E(s)\,ds\right). 
\end{equation}
For $E_3^{c-}$ there is an analogous inequality. The proof relies on the parabolic equation for $D_+D_-\divg(v+2b)$. 

\begin{lemma}
Let $W=W_+=D_+D_-\divg(v+2b)$. Then it holds that
\begin{equation}\label{eq:parW}
\begin{aligned}
D_tW-\calt\Delta W=D_+D_-[D_+,\Delta]\calt&+D_+[D_-,\calt\Delta]\divg(v+2b)\\
&+[D_+,\calt\Delta]D_-\divg(v+2b)+R_w,
\end{aligned}
\end{equation}
where
\begin{equation}\label{eq:Rw}
\begin{aligned}
R_w=&[D_t,D_+]D_-\divg(v+2b)+D_+[D_t,D_-]\divg(v+2b)\\
&+D_+D_-(\divg v\divg(v+b)+2\nabla\calt\cdot\nabla\divg(v+2b)). 
\end{aligned}
\end{equation}
There is an analogous equation for $W_-=D_-D_+\divg(v-2b)$. 
\end{lemma}
\textit{Proof: }It is derived from applying $D_+D_-$ to (\ref{eq:Dtdivvb}) and commuting $D_+D_-$ with $D_t-\calt\Delta$. $\Box$\\[2mm]
The key to proving (\ref{eq:intE3c}) is the following

\begin{proposition}\label{prop:parWform}
It holds that
\begin{equation}\label{eq:parWform}
\calt^{-1}D_tW(D_tW-\calt\Delta W)=\mathcal G+\divg\mathcal F+D_t\mathcal H,
\end{equation}
where $\mathcal G,\mathcal F,\mathcal H$ are subject to the conditions
\begin{equation}\label{eq:Gcond}
\int_\Omega|\mathcal G|\,d\mu_g\leq C(\mathcal K)\left(\mathcal E\log\mathcal E+\sqrt{\mathcal E\log\mathcal E}\Vert D_tW\Vert_{L^2(\Omega)}\right),
\end{equation}
\begin{equation}\label{eq:Fcond}
\mathcal F\cdot N|_{\partial\Omega}=0
\end{equation}
and
\begin{equation}\label{eq:Hcond}
\int_\Omega|\mathcal H|\,d\mu_g\leq \delta\mathcal E(t)+C(\delta^{-1},\mathcal K).
\end{equation}
There is an analogous formulation for $W_-$. 
\end{proposition}
Once we obtain (\ref{eq:parWform})-(\ref{eq:Hcond}), it is rather straightforward to prove (\ref{eq:intE3c}). The rest of this subsection is devoted to validating these conditions. It turns out that most terms in $D_tW-\calt\Delta W$ is bounded in $L^2(\Omega)$ by $\sqrt{\mathcal E\log\mathcal E}$, hence they times $\calt^{-1}D_tW$ can be thrown into $\mathcal G$. This will be made clear in Lemmas \ref{lem:graddiv}-\ref{lem:Rw}. However, there are other terms without such bound, which will be treated in Lemma \ref{lem:PP0}. 

The toughest term in $D_tW-\calt\Delta W$ is $D_+D_-[D_+,\Delta]\calt$. It is decomposed into
\[\begin{aligned}
D_+D_-[D_+,\Delta]\calt=-D_+D_-(\nabla\divg(v+b)\cdot\nabla\calt)+D_+D_-(\curl\curl(v+b)\cdot\nabla\calt)\\
-2D_+D_-(\nabla(v+b)^k\cdot\nabla_k\nabla\calt).
\end{aligned}\]

\begin{lemma}\label{lem:graddiv}
\[
\Vert D_+D_-(\nabla\divg(v+b)\cdot\nabla\calt)\Vert_{L^2(\Omega)}^2\leq C(\mathcal K)\mathcal E.
\]
\end{lemma}
\textit{Proof: }Distributing $D_+D_-$ onto the two factors gives four terms: 
\[\begin{aligned}
D_+D_-\nabla\divg(v+b)\cdot\nabla\calt&+\nabla\divg(v+b)\cdot D_+D_-\nabla\calt\\
&+D_+\nabla\divg(v+b)\cdot D_-\nabla\calt+D_-\nabla\divg(v+b)\cdot D_+\nabla\calt. 
\end{aligned}\] 
The second line is harmless. For the first line, it is important to note the formula
\[
[\nabla,D_+D_-]=\nabla(v+b)^k\nabla_kD_-+\nabla(v-b)^kD_+\nabla_k+D_+\nabla(v-b)^k\nabla_k
\]
together with the fact $\Vert D_+\nabla(v-b)\Vert_{L^p(\Omega)}\leq C(\mathcal K)$, $p<\infty$. Then we may well commute $D_+D_-$ with $\nabla$ in the first line. From a linear combination of (\ref{eq:D3divvb}), (\ref{eq:VpmW2p}) we see that $\Vert\nabla D_+D_-\divg(v+2b)\Vert_{L^2(\Omega)}\leq C(\mathcal K)\sqrt\mathcal E$. Moreover, interpolating between (\ref{eq:DtdivvbLp}) and (\ref{eq:VpmW2p}) gives
\begin{equation}\label{eq:gradVpminterp}
\Vert\nabla D_-\divg(v+2b)\Vert_{L^p(\Omega)}+\Vert\nabla D_+\divg(v-2b)\Vert_{L^p(\Omega)}\leq C(\mathcal K)\mathcal E^\frac{1}{4},\quad p<12.
\end{equation}
Then the result is clear. $\Box$

\begin{lemma}\label{lem:curlcurl}
\[
\Vert D_+D_-(\curl\curl(v+b)\cdot\nabla\calt)-P\Vert_{L^2(\Omega)}^2\leq C(\mathcal K)\mathcal E,
\]
where
\begin{equation}\label{eq:Pdef}
\begin{aligned}
P=D_+\curl\curl(v+b)\cdot D_-\nabla\calt-\nabla(v-b)^i\times D_+\nabla_i\curl(v+b)\cdot\nabla\calt\\
-D_+\curl(\nabla(v-b)^i\times\nabla_i(v+b))\cdot\nabla\calt. 
\end{aligned}
\end{equation}
\end{lemma}
Before proving it, let us apply $\curl$ to (\ref{eq:Dmpvb}) and find that
\begin{equation}\label{eq:curlDpm}
\begin{aligned}
&D_-\curl(v+b)=-\nabla\calt\times\nabla q-\frac{1}{2}\curl[b\divg(v+2b)]-\nabla(v-b)^i\times\nabla_i(v+b),\\
&D_+\curl(v-b)=-\nabla\calt\times\nabla q+\frac{1}{2}\curl[b\divg(v-2b)]-\nabla(v+b)^i\times\nabla_i(v-b).
\end{aligned}
\end{equation}
One may consult Appendix \ref{sec:vectid} for vector identities. \\[2mm]
\textit{Proof of Lemma \ref{lem:curlcurl}: }If $D_+,D_-$ both fall on $\nabla\calt$, the resultant term is harmless, because $\Vert D_+D_-\nabla\calt\Vert_{L^p(\Omega)}\leq C(\mathcal K)\mathcal E^\frac{1}{4}$, $p<12$ due to (\ref{eq:gradVpminterp}), and $\Vert\curl\curl(v+b)\Vert_{L^4(\Omega)}\leq C(\mathcal K)\mathcal E^\frac{1}{4}$. If $D_+$ falls on $\nabla\calt$ and $D_-$ falls on $\curl\curl(v+b)$, the resultant term is also harmless due to (\ref{eq:curlDpm}). There remain two terms: 
\[
D_+\curl\curl(v+b)\cdot D_-\nabla\calt+D_+D_-\curl\curl(v+b)\cdot\nabla\calt.
\]
The first is thrown into $P$. For the second, commute $D_-$ with $\curl$ using $[D_-,\curl]=[D_-,\nabla\times]=-\nabla(v-b)^i\times\nabla_i$. The commutator reads as
\[
-D_+(\nabla(v-b)^i\times\nabla_i\curl(v+b))\cdot\nabla\calt. 
\]
If $D_+$ falls on $\nabla_i\curl(v+b)$, the resultant term is thrown into $P$. If $D_+$ falls on $\nabla(v-b)^i$, then by (\ref{eq:Dmpvb}), (\ref{eq:qW2p}), (\ref{eq:divvbW1p}) the the resultant term is harmless. Finally, we deal with 
\[
D_+\curl D_-\curl(v+b)\cdot\nabla\calt.
\]
Plug (\ref{eq:curlDpm}) into the expression. One readily checks that the only harmful term is $-D_+\curl(\nabla(v-b)^i\times\nabla_i(v+b))\cdot\nabla\calt$, which is thrown into $P$. $\Box$

\begin{lemma}\label{lem:D2T}
\[
\Vert D_+D_-(\nabla(v+b)^k\cdot\nabla_k\nabla\calt)\Vert_{L^2(\Omega)}^2\leq C(\mathcal K)\mathcal E\log\mathcal E. 
\]
\end{lemma}
\textit{Proof: }If $D_-$ falls on $\nabla(v+b)^k$, then by (\ref{eq:Dmpvb}) the resultant term is harmless. The remaining terms are
\[
D_+\nabla(v+b)^k\cdot D_-\nabla_k\nabla\calt+\nabla(v+b)^k\cdot D_+D_-\nabla_k\nabla\calt:=T_1+T_2. 
\]
For $T_1$, let us commute $D_+$ with $\nabla$ and $D_-$ with $\nabla_k$, getting $T_1=\nabla D_+(v+b)^k\cdot\nabla_kD_-\nabla\calt+R$, where $\Vert R\Vert_{L^2(\Omega)}^2\leq C(\mathcal K)\mathcal E$. By interpolation, we have
\[
\Vert\nabla D_+(v+b)^k\Vert_{L^4(\Omega)}^2\leq C(\mathcal K)\Vert D_+(v+b)^k\Vert_{L^\infty(\Omega)}\Vert D_+(v+b)^k\Vert_{H^2(\Omega)}\leq C(\mathcal K)\sqrt\mathcal E\log\mathcal E. 
\]
Then it suffices to show that
\begin{equation}\label{eq:Dmgrad2T}
\Vert\nabla_k D_-\nabla\calt\Vert_{L^4(\Omega)}^2\le \sqrt\mathcal E. 
\end{equation}
This can be a little tricky. Write $D_-\nabla\calt=w+u+\nabla D_-\calt$ where
\[
\left\{\begin{aligned}
&\Delta w=\nabla\divg[D_-,\nabla]\calt-\curl\curl[D_-,\nabla]\calt\quad\text{in }\Omega,\\
&w=0\quad\text{on }\partial\Omega,
\end{aligned}\right.
\]
and
\[
\left\{\begin{aligned}
&\Delta u=0\quad\text{in }\Omega,\\
&u=[D_-,\nabla]\calt\quad\text{on }\partial\Omega. 
\end{aligned}\right.
\]
For $w$, from
\[
-\divg[D_-,\nabla]\calt=\Delta(v-b)^j\partial_j\calt+\nabla(v-b)^j\cdot\nabla\partial_j\calt,
\]
\[
\curl[D_-,\nabla]\calt=\curl(D_-\nabla\calt)=[\curl,D_-]\nabla\calt,
\]
(\ref{eq:Deltadivcurl}) and interpolation it follows that
\[
\Vert\divg[D_-,\nabla]\calt\Vert_{L^4(\Omega)}+\Vert\curl[D_-,\nabla]\calt\Vert_{L^4(\Omega)}\le C(\mathcal K)\mathcal E^{\frac{1}{4}},
\]
then from (\ref{eq:W1pest}) it follows $\Vert\nabla w\Vert_{L^4(\Omega)}\le C(\mathcal K)\mathcal E^{\frac{1}{4}}$. For $u$, from 
\[
\Vert u-[D_-,\nabla]\calt\Vert_{H^2(\Omega)}\le C(\mathcal K)\Vert\Delta[D_-,\nabla]\calt\Vert_{L^2(\Omega)}\le C(\mathcal K)\sqrt\mathcal E,
\]
(\ref{eq:vBH3}), (\ref{eq:TH3}) we see $\Vert u\Vert_{H^2(\Omega)}\le C(\mathcal K)\sqrt\mathcal E$, and from $\Vert u\Vert_{L^\infty(\partial\Omega)}\leq C(\mathcal K)$ (using (\ref{eq:gradvBbdry}), (\ref{eq:TW2p})), the maximal principle we see $\Vert u\Vert_{L^\infty(\Omega)}\leq C(\mathcal K)$. Then by interpolation we have $\Vert\nabla w\Vert_{L^4(\Omega)}\le C(\mathcal K)\mathcal E^{\frac{1}{4}}$, and (\ref{eq:Dmgrad2T}) follows. 

For $T_2=\nabla(v+b)^k\cdot D_+D_-\nabla_k\nabla\calt$, write $D_+D_-\nabla_k\nabla\calt=D_+[D_-,\nabla_k\nabla]\calt+[D_+,\nabla_k\nabla]D_-\calt+\nabla_k\nabla D_+D_-\calt$; it is important that $D_+$ matches $v-b$ in the first term. Then by (\ref{eq:Dmpvb}), $D_-\calt=\calt\divg(v-2b)$ and Proposition \ref{prop:ellipest} we have 
\[
\Vert T_2\Vert_{L^2(\Omega)}\leq C(\mathcal K)(\Vert|\nabla(v-b)|\nabla^3q\Vert_{L^2(\Omega)}+\sqrt{\mathcal E}). 
\]
Next, invoking (\ref{eq:wHr}) with $\rho=|\nabla(v-b)|$ and combining (\ref{eq:gradvBbdry}), (\ref{eq:DeltaqW16}), we see $\Vert|\nabla(v-b)|\nabla^3q\Vert_{L^2(\Omega)}\leq C(\mathcal K)\sqrt{\mathcal E}$. The proof is now completed. $\Box$\\
Leaving $P$ aside, we move on to control the other terms in $D_tW-\calt\Delta W$. 

\begin{lemma}\label{lem:Rw}
Recall (\ref{eq:parW}). For the second line we have
\[
\Vert [D_+,\calt\Delta]D_-\divg(v+2b)\Vert_{L^2(\Omega)}^2+\Vert R_w\Vert_{L^2(\Omega)}^2\leq C(\mathcal K)\mathcal E. 
\]
For the remaining $D_+[D_-,\calt\Delta]\divg(v+2b)$, we have
\[
\Vert D_+[D_-,\calt\Delta]\divg(v+2b)+2P_0\Vert_{L^2(\Omega)}^2\leq C(\mathcal K)\mathcal E\log\mathcal E, 
\]
where
\begin{equation}\label{eq:P0def}
P_0=\calt\divg(\nabla(v-b)^iD_+\nabla_i\divg(v+2b)). 
\end{equation}
\end{lemma}
\textit{Proof: }The bound for $[D_+,\calt\Delta]D_-\divg(v+2b)$ follows (\ref{eq:VpmW2p}), (\ref{eq:gradVpminterp}) and the interpolation estimate of $\Delta(v+b)$ using (\ref{eq:Deltadivcurl}). The estimate of $R_w$ is even easier. For $D_+[D_-,\calt\Delta]\divg(v+2b)$, write
\[\begin{aligned}
[D_-,\calt\Delta]\divg(v+2b)=(D_-\calt)\Delta\divg(v+2b)+\calt\Delta(v-b)^i\nabla_i\divg(v+2b)\\
-2\calt\divg(\nabla(v-b)^i\nabla_i\divg(v+2b)). 
\end{aligned}\]
Apply $D_+$ to the three terms on the right. The first term, $D_+[(D_-\calt)\Delta\divg(v+2b)]$, is harmless due to (\ref{eq:D3divvb}). For the second, distributing $D_+$ to its factors gives 
\[
D_+(\calt\Delta(v-b)^i)\nabla_i\divg(v+2b)+\calt\Delta(v-b)^iD_+\nabla_i\divg(v+2b):=T_1+T_2.
\]
For $T_1$ we use (\ref{eq:Dmpvb}), (\ref{eq:DeltaqW16}). For $T_2$ we use interpolation estimates of $\Delta(v-b)^i$, $\nabla_iD_+\divg(v+2b)$; the $\log\mathcal E$ factor arises here due to
\begin{equation}\label{eq:Dpdivpinterp}
\Vert\nabla_iD_+\divg(v+2b)\Vert_{L^4(\Omega)}^2\leq C(\mathcal K)\sqrt{\mathcal E\log\mathcal E}, 
\end{equation}
as is seen from $D_+=D_-+2b^k\nabla_k$, (\ref{eq:gradVpminterp}) and (\ref{eq:divvbLip}), (\ref{eq:D3divvb}). Finally, we treat $D_+\divg(\nabla(v-b)^i\nabla_i\divg(v+2b))$. Commute derivatives and distribute $D_+$ onto the two factors inside. As is argued above, the commutator and the term containing $D_+\nabla(v-b)^i$ are harmless. The remaining term is put into $P_0$. $\Box$\\[2mm]
Till now we have control of $(D_t-\calt\Delta)W$ except for $P,P_0$ in (\ref{eq:Pdef}), (\ref{eq:P0def}). They will be treated in the following 

\begin{lemma}\label{lem:PP0}
\[
\calt^{-1}D_tW(P-2P_0)=\mathcal G+\divg\mathcal F+D_t\mathcal H,
\]
where $\mathcal F,\mathcal G,\mathcal H$ are as (\ref{eq:Gcond})-(\ref{eq:Hcond}). 
\end{lemma}
\textit{Proof: }We first consider $P_0$. Write $F_0=\nabla(v-b)^iD_+\nabla_i\divg(v+2b)$ so that
\[\begin{aligned}
\calt^{-1}D_tWP_0=&D_tW\divg F_0 \\
=&\divg(D_tWF_0)-D_t(\nabla W\cdot F_0)+\nabla W\cdot D_tF_0-[\nabla,D_t]W\cdot F_0.
\end{aligned}\]
The first term can be attributed to $\divg\mathcal F$ due to $D_+W|_{\partial\Omega}=0$. The last two terms can be attributed to $\mathcal G$ due to (\ref{eq:D3divvb}), (\ref{eq:VpmW2p}), (\ref{eq:Dpdivpinterp}) and interpolation estimate of $\nabla D_t(v-b)^i$; it is also important to note $\Vert\nabla W\Vert_{L^2(\Omega)}^2=E_3^{c+}$. The second term is attributed to $D_t\mathcal H$ due to (\ref{eq:Dpdivpinterp}) and convex inequality. 

We next consider $P$. The proof is quite similar, thus we only describe it briefly. Decompose $P$ into
\[
P_1=D_+\curl\curl(v+b)\cdot D_-\nabla\calt,\quad P_2=D_+\curl(\nabla(v-b)^i\times\nabla_i(v+b))\cdot\nabla\calt,
\]
\[
P_3=-\nabla(v-b)^i\times D_+\nabla_i\curl(v+b)\cdot\nabla\calt=\nabla(v-b)^i\times \nabla\calt\cdot D_+\nabla_i\curl(v+b). 
\]
For $P_1$, we have
\[\begin{aligned}
\calt^{-1}D_tWP_1=&\curl D_+\curl(v+b)\cdot\nabla(\calt^{-1}D_-\calt)\,D_tW+\mathcal G\\
=&\divg(D_+\curl(v+b)\times\nabla(\calt^{-1}D_-\calt))\,D_tW+\mathcal G\nonumber\\
%=&-D_+\curl(v+b)\times\nabla(\calt^{-1}D_-\calt)\cdot D_t\nabla W+\divg\mathcal F+\mathcal G\nonumber\\
=&-D_t(D_+\curl(v+b)\times\nabla(\calt^{-1}D_-\calt)\cdot \nabla W)\\
&+D_t(D_+\curl(v+b)\times\nabla(\calt^{-1}D_-\calt))\cdot\nabla W+\divg\mathcal F+\mathcal G.
%=&D_t\mathcal H+\divg\mathcal F+\mathcal G. 
\end{aligned}\]
Thanks to (\ref{eq:D3vB}), (\ref{eq:divvbLip}) and $\calt^{-1}D_-\calt=\divg(v-2b)$, we eventually get $D_t\mathcal H+\divg\mathcal F+\mathcal G$. 
%In (\ref{eq:P1line1}) we have used that $D_t=D_-+b^k\nabla_k$ and used (\ref{eq:D3vB}), (\ref{eq:divvbLip}) to control $\curl D_+\curl(v+b)\cdot\nabla(\calt^{-1}D_-\calt)\,b^k\nabla_kW$. As discussed for $P_0$, (\ref{eq:P1line4}) has the form $D_t\mathcal H+\divg\mathcal F+\mathcal G$. The first term in (\ref{eq:P1line5}) is controlled like $\mathcal G$. 
For $P_2$, we have
\[\begin{aligned}
\calt^{-1}D_tWP_2=&\curl D_+(\nabla(v-b)^i\times\nabla_i(v+b))\cdot\nabla\log\calt\,D_tW+\mathcal G\\
=&\divg(D_+(\nabla(v-b)^i\times\nabla_i(v+b))\times\nabla\log\calt)\,D_tW+\mathcal G\\
%=&-D_+(\nabla(v-b)^i\times\nabla_i(v+b))\times\nabla\log\calt\cdot D_t\nabla W+\divg\mathcal F+\mathcal G\\
=&-D_t(D_+(\nabla(v-b)^i\times\nabla_i(v+b))\times\nabla\log\calt\cdot \nabla W)\\
&+D_t(D_+(\nabla(v-b)^i\times\nabla_i(v+b))\times\nabla\log\calt)\cdot \nabla W+\divg\mathcal F+\mathcal G,
\end{aligned}\]
which is again $D_t\mathcal H+\divg\mathcal F+\mathcal G$. For $P_3$, noting $D_+=D_-+2b^k\nabla_k$ we have
\[\begin{aligned}
\calt^{-1}D_tWP_2=&\nabla(v-b)^i\times \nabla\log\calt\cdot D_-\nabla_i\curl(v+b)\,D_tW\\
&+2\nabla(v-b)^i\times \nabla\calt\cdot b^k\nabla_k\nabla_i\curl(v+b)\,D_tW. 
\end{aligned}\]
Thanks to (\ref{eq:curlDpm}), the first line can be attributed to $\mathcal G$. Next, noting $(b^k\nabla_kf)h=\divg(fhb)-fb^k\nabla_kh-fh\divg b$ and $b\cdot N|_{\partial\Omega}=0$, we have
%then we continue to write
\[\begin{aligned}
\calt^{-1}D_tWP_2=&-2b^k\nabla_k(\nabla(v-b)^i\times \nabla\calt)\cdot\nabla_i\curl(v+b)\,D_tW\\
&-2\nabla(v-b)^i\times \nabla\calt\cdot\nabla_i\curl(v+b)\,D_t(b^k\nabla_kW)+\divg\mathcal F+\mathcal G. 
\end{aligned}\]
The first line can be thrown into $\mathcal G$ due to $\Vert\nabla^2(v-b)^i\Vert_{L^4(\Omega)}^2\leq C(\mathcal K)\sqrt{\mathcal E}\log\mathcal E$ and $\Vert\nabla_i\curl(v+b)\Vert_{L^4(\Omega)}^2\leq C(\mathcal K)\sqrt{\mathcal E}$. For the second line, using $fD_th=D_t(fh)-(D_tf)h$ we have
\[\begin{aligned}
\calt^{-1}D_tWP_2=&D_t(-2\nabla(v-b)^i\times \nabla\calt\cdot\nabla_i\curl(v+b)\,b^k\nabla_kW)\\
&+D_t(2\nabla(v-b)^i\times \nabla\calt\cdot\nabla_i\curl(v+b))b^k\nabla_kW+\divg\mathcal F+\mathcal G, 
\end{aligned}\]
which is again $D_t\mathcal H+\divg\mathcal F+\mathcal G$. This completes the proof. $\Box$\\[2mm]
\textit{Proof of (\ref{eq:intE3c}): }Proposition \ref{prop:parWform} follows from the preceding lemmas. Based on this, integrating (\ref{eq:parWform}) over $\Omega$ and using the divergence theorem, we get
\[\begin{aligned}
\frac{d}{dt}\frac{1}{2}\int_\Omega|\nabla W|^2\,d\mu_g=&\frac{d}{dt}\int_\Omega\mathcal H\,d\mu_g-\int_\Omega\calt^{-1}(D_tW)^2\,d\mu_g\\
&+\int_\Omega\mathcal G-[\nabla,D_t]W\cdot\nabla W+(\frac{1}{2}|\nabla W|^2-\mathcal H)\divg v\,d\mu_g. 
\end{aligned}\]
Further integrating in $t$ and using the bounds (\ref{eq:Gcond}), (\ref{eq:Hcond}), we obtain (\ref{eq:intE3c}). $\Box$

\subsection{Estimates of temporal derivatives of $E_3^a,E_3^b$}

We shall show that
\begin{equation}\label{eq:dtE3a}
\frac{d}{dt}E_3^a(t)\leq C(\mathcal K)\mathcal E(\Vert\nabla_ND_tq\Vert_{L^\infty(\partial\Omega)}+\log\mathcal E),
\end{equation}
\begin{equation}\label{eq:dtE3b}
\frac{d}{dt}E_3^b(t)\leq C(\mathcal K)\mathcal E\log\mathcal E. 
\end{equation}
\textit{Proof of (\ref{eq:dtE3a}): }Take $r=3$ in (\ref{eq:dtEra}) and recall (\ref{eq:I1adef})-(\ref{eq:I3adef}). By (\ref{eq:divvbH3}), (\ref{eq:D3vB}) and interpolation we bound $\mathcal I_3^a$ together with the first line in (\ref{eq:dtEra}) by $C(\mathcal K)\mathcal E\log\mathcal E$. For $\mathcal I_1^a$, $\mathcal I_1^a$, note from Lemma \ref{lem:Nmu} (together with the remark below), the calculation $D_t\mathfrak a=\mathfrak a^2(\nabla_ND_tq-\nabla_Nv\cdot\nabla q)$ and the estimates (\ref{eq:vBLip}), (\ref{eq:gradvBbdry}) that
\[
|Q_t|+|Q_{,k}|\leq C(\mathcal K)\log\mathcal E\quad\text{in }\Omega,
\]\[
\left|\frac{D_t(\mathfrak ad\mu_\zeta)}{\mathfrak ad\mu_\zeta}\right|\leq C(\mathcal K)(M+\Vert\nabla_ND_tq\Vert_{L^\infty(\partial\Omega)})\quad\text{on }\partial\Omega. 
\]
Further note from (\ref{eq:QPi}) that $Q_t(\nabla^rq,\nabla^rq)=2\langle\Pi_t(\nabla^rq),\Pi(\nabla^rq)\rangle$. Then from (\ref{eq:D3q}), (\ref{eq:D2qNbdry}), (\ref{eq:D3qbdry}) we see that $|\mathcal I_1^a|+|\mathcal I_2^a|\leq C(\mathcal K)\mathcal E(\Vert\nabla_ND_tq\Vert_{L^\infty(\partial\Omega)}+\log\mathcal E)$. 

We finally control $Q(D_t\nabla^3q-\mathfrak a^{-1}N^k\nabla^3v_k,\nabla^3q)$. In view of (\ref{eq:D2qNbdry}) and (\ref{eq:hypq}), (\ref{eq:PiD3Dtqbdry}), it remains to bound the tangential projection of 
\[
A=\sum_{s=1}^{2}\binom{3}{s}\nabla^sv^k\widetilde\otimes\nabla^{3-s}\nabla_kq
\]
in $L^2(\partial\Omega)$. By H\"older inequality, we have
\[\begin{aligned}
\Vert A\Vert_{L^2(\partial\Omega)}\leq&\Vert\nabla v\Vert_{L^\infty(\partial\Omega)}\Vert\nabla^3q\Vert_{L^2(\partial\Omega)}+\Vert\nabla^2v\Vert_{L^3(\partial\Omega)}\Vert\nabla^2q\Vert_{L^6(\partial\Omega)}.
\end{aligned}\]
For the first term we use (\ref{eq:gradvBbdry}), (\ref{eq:D3qbdry}). For the second term we use (\ref{eq:gradinterp}), (\ref{eq:trinterp}) and (\ref{eq:vBLip}), (\ref{eq:qW2p}), (\ref{eq:vBH3}), (\ref{eq:D3q}). Now (\ref{eq:dtE3a}) is clear. $\Box$\\[2mm]
\textit{Proof of (\ref{eq:dtE3b}): }Take $r=3$ in (\ref{eq:dtErb}) and recall (\ref{eq:Ibdef}). By (\ref{eq:curlv}), (\ref{eq:curlB}), Propositions \ref{prop:parest}, \ref{prop:ellipest} and (\ref{eq:gradinterp}), we have
\[
\Vert\calt^{-1}D_t\curl v-B^i\nabla_i\curl B\Vert_{H^2(\mathcal D_t)}^2\le C(\mathcal K)\mathcal E\log\mathcal E,
\]
\[
\Vert D_t\curl B-B^i\nabla_i\curl v\Vert_{H^2(\mathcal D_t)}\le C(\mathcal K)\mathcal E\log\mathcal E,
\]
and
\[
|\mathcal I^b|\le C(\mathcal K)\mathcal E\log\mathcal E. 
\]
This completes the proof. $\Box$

\subsection{Completion of the proof}

In the preceding subsections we control the growth of $E_3$. For $E_1,E_2$, one readily checks that $\frac{d}{dt}(E_1(t)+E_2(t))\leq C(\mathcal K)\mathcal E$. Denote $\rho(t)=\sup_{s\leq t}\mathcal E(s)$ and $\gamma(t)=1+\Vert\nabla_ND_tq\Vert_{L^\infty(\partial\Omega)}$, we then have
\[
\rho(t)\leq C(\mathcal K)\left(\rho(0)+\int_0^t\gamma(s)\log(\rho(s))\rho(s)\,ds\right). 
\]
Moreover, when $B=0$ we may set $\gamma=1$ thanks to (\ref{eq:gradDtqap}). Since $\mu(r)=r\log\frac{1}{r},\ r\in[0,1]$ is an Osgood modulus of continuity (Definition 3.1 in  \cite{FAPDE}), from the dual Osgood lemma (Lemma 3.8 in \cite{FAPDE}) it follows that $\mathcal E(t)\le C(\mathcal K,\overline L,\mathcal E(0))$, $t\in[0,T)$. Recall the definition of $\mathcal E_{r_*}$ from (\ref{eq:Er*def}) with $r_*=[\frac{n}{2}]+2=3$, then by Proposition \ref{prop:ellipest} we also have 
\[
\Vert\nabla^2q\Vert_{L^\infty(\Omega)}^2+\Vert\nabla_ND_tq\Vert_{L^\infty(\partial\Omega)}^2+\mathcal E_{r_*}(t)\le C(\mathcal K,\overline L,\mathcal E(0)),\quad t\in[0,T). 
\]
As mentioned above, when $B=0$ there is no $\overline L$ on the right. Now Theorem \ref{thm:buc} follows from (\ref{eq:dtErest}) combined with an induction on $r>r_*$.

\section{Regularity estimates II}

In this section we derive estimates needed for proving Theorem \ref{thm:buc}. In section \ref{sec:parest} basic controls are derived from the parameter $\mathcal K$, and in section \ref{sec:ellipest} more regularity is recovered from the energy $\mathcal E$. 

\subsection{Parabolic estimates}\label{sec:parest}

In this subsection we derive the following estimates. 

\begin{proposition}\label{prop:parest}
(i) For $\vol\Omega$ and $\calt$, we have
\begin{equation}\label{eq:TTinv}
C(L_0)^{-1}\leq\sup_\Omega\left(\calt+\frac{1}{\calt}\right)\leq C(L_0),
\end{equation}
\begin{equation}\label{eq:volume}
C(L_0)^{-1}\leq\vol\Omega\leq C(L_0),
\end{equation}
and
\begin{equation}\label{eq:TW2p}
\Vert\nabla\calt\Vert_{L^\infty(\Omega)}+\Vert\nabla^2\calt\Vert_{L^p(\Omega)}\leq C(K,L_0)M,\quad p<\infty.
\end{equation}
Here $L_0$ is a tuple of initial quantities: 
\begin{equation}\label{eq:L0def}
L_0=\left.\left(\inf_{\Omega}\calt,\ \sup_{\Omega}\calt,\ \int_{\Omega}\calt^{-1}d\mu_g,\ \frac{1}{\vol\Omega}\int_{\Omega}Bd\mu_g\right)\right|_{t=0}.
\end{equation}
(ii) For $v,B$ we have
\begin{equation}\label{eq:BLinf}
\Vert B\Vert_{L^\infty(\Omega)}\leq C(K,L_0)M,
\end{equation}
\begin{equation}\label{eq:vBW1p}
\Vert\nabla v\Vert_{L^p(\Omega)}+\Vert\nabla B\Vert_{L^p(\Omega)}\leq C(K,L_0)M,\quad p<\infty,
\end{equation}
and
\begin{equation}\label{eq:vBLip}
\Vert\nabla v\Vert_{L^\infty(\Omega)}+\Vert\nabla B\Vert_{L^\infty(\Omega)}\leq C(K,M,L_0)\log\mathcal E. 
\end{equation}
\begin{equation}\label{eq:gradvBbdry}
\Vert\nabla v\Vert_{L^\infty(\partial\Omega)}+\Vert\nabla B\Vert_{L^\infty(\partial\Omega)}\leq C(K,M,L_0). 
\end{equation}
(iii) Recall $b=\calt^\frac{1}{2}B$. For $\divg v,\divg b$ we have
\begin{equation}\label{eq:DtdivvbLp}
\Vert D_t\divg v\Vert_{L^p(\Omega)}+\Vert D_t\divg b\Vert_{L^p(\Omega)}\leq C(K,M,L_0,L_1),\quad p<\infty,
\end{equation}
\begin{equation}\label{eq:divvbW1p}
\Vert\nabla\divg v\Vert_{L^p(\Omega)}+\Vert\nabla\divg b\Vert_{L^p(\Omega)}\leq C(K,M,L_0,L_1),\quad p<\infty,
\end{equation}
and
\begin{equation}\label{eq:divvbLip}
\Vert\nabla\divg v\Vert_{L^\infty(\Omega)}+\Vert\nabla\divg b\Vert_{L^\infty(\Omega)}\leq C(K,M,L_0,L_1)\log\mathcal E.
\end{equation}
Here $L_1$ is the initial quantity defined by
\begin{equation}\label{eq:L1def}
L_1=(\Vert D_-\divg(v+2b)\Vert_{L^\infty(\Omega)}+\Vert D_+\divg(v-2b)\Vert_{L^\infty(\Omega)})|_{t=0}.
\end{equation}
(iv) For the total pressure $q$ we have
\begin{equation}\label{eq:qW2p}
\Vert\nabla q\Vert_{L^\infty(\Omega)}+\Vert\nabla^2q\Vert_{L^p(\Omega)}\leq C(K,M,L_0,L_1),\quad p<\infty.
\end{equation}
(v) When $B=0$, we additionally have
\begin{equation}\label{eq:Dt2divvap}
\Vert D_t^2\divg v\Vert_{L^4(\Omega)}\leq C(K,M,L_0,L_1,\Vert D_t^2\divg v\Vert_{L^4(\Omega)}|_{t=0})
\end{equation}
and
\begin{equation}\label{eq:gradDtqap}
\Vert\nabla D_tq\Vert_{L^\infty(\Omega)}\leq C(K,M,L_0,L_1,\Vert D_t^2\divg v\Vert_{L^4(\Omega)}|_{t=0})\log\mathcal E.
\end{equation}
\end{proposition}
\textit{Proof of (\ref{eq:TTinv})-(\ref{eq:volume}): }The estimate of $\calt,\calt^{-1}$ follows from the maximal principle for the parabolic equation $D_t\calt=\calt\Delta\calt$ with Dirichlet boundary condition $\calt|_{\partial\Omega}=1$. For (\ref{eq:volume}), note the conservation law
\[
\frac{d}{dt}\int_\Omega\calt^{-1}d\mu_g=\int_\Omega D_t(\calt^{-1})+\divg v\,d\mu_g=0,
\]
then the result follows from (\ref{eq:TTinv}). $\Box$\\[2mm]
\textit{Proof of (\ref{eq:TW2p})-(\ref{eq:gradvBbdry}): }Inequality (\ref{eq:TW2p}) follows from $\Delta\calt=\divg v$ and (\ref{eq:W2pest}). For (\ref{eq:vBW1p}), note the equation $\Delta v=\nabla\divg v-\curl\curl v$, then the result follows from (\ref{eq:W1pest}), (\ref{eq:gradLp}) by an obvious decomposition of the elliptic problem into harmonic part and homogeneous boundary condition part. For (\ref{eq:vBLip}), we use (\ref{eq:gradest}), (\ref{eq:gradest2}) and the regularity bound (\ref{eq:vBH3}). For (\ref{eq:gradvBbdry}), we use (\ref{eq:ptHodge}). For (\ref{eq:BLinf}), we first note from $\divg B=0$ and $B_kN^k|_{\partial\Omega}=0$ that
\[ 
\frac{d}{dt}\int_\Omega B\,d\mu_g=\int_\Omega D_tB+B\divg v\,d\mu_g=\int_\Omega B^k\nabla_kv\,d\mu_g=0.
\]
Then it remains to control $\text{osc}_\Omega B$, for which we use Morrey inequality. However, we have to mind the Sobolev constant here. For this, we switch to Eulerian coordinates and select $m$ balls $B_\frac{1}{K}(x_i)$ such that $\mathcal D_t\subset\cup B_\frac{1}{K}(x_i)$ and that $B_\frac{1}{3K}(x_i)$ are disjoint. This implies $m\leq CK^3\vol\mathcal D_t$. In each $B_\frac{1}{K}(x_i)\cap\mathcal D_t$, the local geometry guarantees that any pair of points $x,x'\in B_r(x_i)\cap\mathcal D_t$ can be joint by a curve contained in $\mathcal D_t$ of length less or equal than $4|x-x'|$. This enable us to use the following version of Morrey inequality \cite{AnaXd}:
\[
|B(t,x)-B(t,x')|\leq C|x-x'|^{\frac{1}{2}}\Vert\nabla B\Vert_{L^6(B_r(x_i)\cap\mathcal D_t)},\quad x,x'\in B_r(x_i)\cap\mathcal D_t, 
\] 
from which (\ref{eq:BLinf}) follows. $\Box$\\[2mm]
\textbf{Remark. }There is no danger of circular argument in the proof of (\ref{eq:vBLip}), because we do not use it again for deriving (\ref{eq:vBH3}). \\[2mm]
The proof of (\ref{eq:DtdivvbLp}) relies on the parabolic equations (\ref{eq:Dtdivvb}), (\ref{eq:Dtdivvb1}), but in divergence form, as following. 

\begin{lemma}
(i) 
\begin{equation}\label{eq:divTdq}
-\divg(\calt\nabla q)=D_t\divg v+\partial_iv^j\partial_jv^i-\calt\partial_iB^j\partial_jB^i-\partial_i\calt B^j\partial_jB^i. 
\end{equation}
(ii) 
\begin{equation}\label{eq:pardivvb}
\begin{aligned}
D_t\divg(v+2b)-\divg(\calt\nabla\divg(v+2b))=\divg G_+ -(\divg v)\divg(v+2b),\\
D_t\divg(v-2b)-\divg(\calt\nabla\divg(v-2b))=\divg G_- -(\divg v)\divg(v-2b),
\end{aligned}
\end{equation}
where
\[
G_+=\divg(2v+3b)\nabla\calt-\sum_{j=1}^3\nabla_j\calt(\nabla_j(v+b)+\nabla(v+b)_j),
\]
\[
G_-=\divg(2v-3b)\nabla\calt-\sum_{j=1}^3\nabla_j\calt(\nabla_j(v-b)+\nabla(v-b)_j).
\]
(iii) Let $V_+=D_-\divg(v+2b)$, $V_-=D_+\divg(v-2b)$. Then we have
\begin{equation}\label{eq:parVpm}
\begin{aligned}
D_tV_+-\divg(\calt\nabla V_+)=D_-\divg G_++\divg F_++R_+,\\
D_tV_--\divg(\calt\nabla V_-)=D_+\divg G_-+\divg F_-+R_-,
\end{aligned}
\end{equation}
where
\[
F_+=[D_-,\calt\nabla]\divg(v+2b)-\sum_{i=1}^3\calt\nabla_i\divg(v+2b)\nabla_i(v-b)+\divg(v-b)\calt\nabla\divg(v+2b),
\]
\[
F_-=[D_+,\calt\nabla]\divg(v-2b)-\sum_{i=1}^3\calt\nabla_i\divg(v-2b)\nabla_i(v+b)+\divg(v+b)\calt\nabla\divg(v-2b),
\]
and
\[
R_+=[D_t,D_-]\divg(v+2b)-D_-(\divg v\divg(v+2b))-\divg(v-b)\divg(\calt\nabla\divg(v+2b)),
\]
\[
R_-=[D_t,D_+]\divg(v-2b)-D_+(\divg v\divg(v-2b))-\divg(v+b)\divg(\calt\nabla\divg(v-2b)).
\]
\end{lemma}
\textit{Proof: }(i) This is just (\ref{eq:ellipdtsq2}) with $s=0$. \\
(ii) We only prove (\ref{eq:pardivvb}$+$). Recall (\ref{eq:parmid}):
\[
(\divg v)\divg b+D_t\divg(v+2b)=D_+\Delta\calt. 
\]
The right-hand side equals
\[
[D_+,\divg]\nabla\calt+\divg([D_+,\nabla]\calt)+\divg(\divg(v+2b)\nabla\calt)+\divg(\calt\nabla\divg(v+2b)).
\]
Using $[D_+,\nabla]=-\nabla(v+b)^k\nabla_k$,
\begin{equation}\label{eq:[Dp,div]}
[D_+,\divg]F=-\nabla_j(F\cdot\nabla(v+b)^j)+\divg(\divg(v+b)F)-\divg(v+b)\divg F
\end{equation}
for $F=\nabla\calt$ and $\Delta\calt=\divg v$, we obtain the desired identity. \\
(iii) Still, we only prove (\ref{eq:parVpm}$+$). Apply $D_-$ to (\ref{eq:pardivvb}$+$) and commute derivatives, we see
\[\begin{aligned}
D_tV_+-\divg(\calt\nabla V_+)=D_-\divg G_+-D_-(\divg v\divg(v+2b))+[D_t,D_-]\divg(v+2b)\\
+[D_-,\divg](\calt\nabla\divg(v+2b))+\divg([D_-,\calt\nabla]\divg(v+2b)). 
\end{aligned}\]
Again, using (\ref{eq:[Dp,div]}) with $F=\calt\nabla\divg(v+2b)$ and $D_+$ replaced by $D_-$ and correspondingly, $v+b$ replaced by $v-b$, we obtain (\ref{eq:parVpm}$+$). $\Box$

\begin{lemma}\label{lem:pardiv}
Suppose $V$ solves
\begin{equation}\label{eq:pardiv}
D_tV-\divg\left(e^f\nabla V\right)=\divg\mathcal F+\mathcal R\quad\text{in }[0,T)\times\Omega
\end{equation}
subject to $V|_{\partial\Omega}=0$ and let $p$ be a positive even integer. Suppose there is a constant $\mathcal A>0$ such that
\begin{equation}\label{eq:fdivvcond}
\Vert f\Vert_{L^\infty(\Omega)}+\Vert \divg v\Vert_{L^\infty(\Omega)}\leq\mathcal A
\end{equation}
and
\begin{equation}\label{eq:FRcond}
\Vert\mathcal F\Vert_{L^p(\Omega)}+\Vert\mathcal R\Vert_{L^p(\Omega)}\leq\mathcal A\Vert V\Vert_{L^p(\Omega)}.
\end{equation}
Then we have
\[
\sup_{t\in[0,T)}\Vert V\Vert_{L^p(\Omega)}\leq C(T,\mathcal A)(\Vert V\Vert_{L^p(\Omega)}|_{t=0}+T). 
\]
\end{lemma}
\textit{Proof: }Multiplicating both sides in (\ref{eq:pardiv}) by $V^{p-1}$ and integrating by parts, we get
\[\begin{aligned}
&\frac{d}{dt}\frac{1}{p}\int_\Omega V^p\,d\mu_g+\frac{4(p-1)}{p^2}\int_\Omega e^f|\nabla(V^\frac{p}{2})|^2\,d\mu_g\\
&=\int_\Omega(\mathcal R-\frac{1}{p}\divg v)V^{p-1}\,d\mu_g-\frac{2(p-1)}{p}\int_\Omega\mathcal F\cdot\nabla(V^\frac{p}{2})V^{\frac{p}{2}-1}\,d\mu_g. 
\end{aligned}\]
Then the result follows from H\"older and Gronwall inequalities. $\Box$\\[2mm]
\textit{Proof of (\ref{eq:DtdivvbLp})-(\ref{eq:divvbLip}), (\ref{eq:qW2p}): }We first control $V_+,V_-$ in $L^p(\Omega)$ for any even integer $p\geq4$. In view of the above lemma, we only need show that (\ref{eq:parVpm}) has the form (\ref{eq:pardiv})-(\ref{eq:FRcond}). Condition (\ref{eq:fdivvcond}) is clearly satisfied with $e^f=\calt$ and $\mathcal A=(K,M,L_0)$. For (\ref{eq:FRcond}), we need some elliptic bounds: 
\begin{equation}\label{eq:divvbW1ptd}
\Vert\nabla\divg v\Vert_{L^{\tilde p}(\Omega)}+\Vert\nabla\divg b\Vert_{L^{\tilde p}(\Omega)}\leq C(\mathcal A)(1+\Vert V_+\Vert_{L^4(\Omega)}+\Vert V_-\Vert_{L^4(\Omega)}),\quad \tilde p<\infty, 
\end{equation}
\begin{equation}\label{eq:qW2ptemp}
\Vert\nabla q\Vert_{L^\infty(\Omega)}+\Vert\nabla^2 q\Vert_{L^p(\Omega)}\leq C(\mathcal A)(1+\Vert V_+\Vert_{L^p(\Omega)}+\Vert V_-\Vert_{L^p(\Omega)}). 
\end{equation}
To prove (\ref{eq:divvbW1ptd}), we use (\ref{eq:pardivvb}). Add both sides by $-b^k\nabla_k\divg(v+2b)$, note $D_t-b^k\nabla_k=D_-$, $b^k\nabla_k\divg(v+2b)=\divg(\divg(v+2b)b)-\divg(v+2b)\divg b$ and invoke (\ref{eq:W1pest}), we then obtain $\Vert\nabla\divg(v\pm2b)\Vert_{L^{\tilde p}(\Omega)}\leq C(\mathcal A)(1+\Vert V_\pm\Vert_{L^p(\Omega)})$, and (\ref{eq:divvbW1ptd}) follows from a linear combination. As a consequence, from the formulas (\ref{eq:Dtdivvrec}), (\ref{eq:Dtdivbrec}): 
\[%\begin{equation}\label{eq:Dtdivvrec}
D_-\divg(v+2b)+D_+\divg(v-2b)=2D_t\divg v-4b^k\nabla_k\divg b,
\]%\end{equation}
\[%\begin{equation}\label{eq:Dtdivbrec}
D_-\divg(v+2b)-D_+\divg(v-2b)=4D_t\divg b-2b^k\nabla_k\divg v,
\]%\end{equation}
it follows that
\begin{equation}\label{eq:DtdivvbLptemp}
\Vert D_t\divg v\Vert_{L^p(\Omega)}+\Vert D_t\divg b\Vert_{L^p(\Omega)}\leq C(\mathcal A)(1+\Vert V_+\Vert_{L^p(\Omega)}+\Vert V_-\Vert_{L^p(\Omega)}).
\end{equation}
Then (\ref{eq:qW2ptemp}) follows from (\ref{eq:divTdq}) and (\ref{eq:W2pest}). 

Now we are in a position to establish (\ref{eq:FRcond}) for (\ref{eq:parVpm}); as before, we only consider the ($+$) part. Clearly, $F_+$ can be attributed to $\mathcal F$ and the first two terms in $R_+$ can be attributed to $\mathcal R$. One may worry about the third term, but from (\ref{eq:pardivvb}) it follows that $\divg(\calt\nabla\divg(v+2b))$ has the estimate as $\mathcal R$ as well. It remains to control $D_-\divg G_+$. According to the formula
\begin{equation}\label{eq:commdivform}
[\nabla_i,D_\pm]\phi=\nabla_j(\phi\nabla_i(v\pm b)^j)-\phi\nabla_i\divg(v\pm b),
\end{equation}
the commutator $[D_-,\divg]G_+$ is controlled as $\divg\mathcal F+\mathcal R$. Thanks to (\ref{eq:DtdivvbLptemp}), (\ref{eq:qW2ptemp}) and (\ref{eq:Dmpvb}), the vector $D_-G_+$ is controlled as $\mathcal F$. Then it follows from Lemma \ref{lem:pardiv} that
\[
\Vert V_+\Vert_{L^p(\Omega)}+\Vert V_-\Vert_{L^p(\Omega)}\leq C(K,M,L_0,T)(L_1+T). 
\]
Inequalities (\ref{eq:DtdivvbLp})-(\ref{eq:divvbW1p}), (\ref{eq:qW2p}) follow from the above one combined with (\ref{eq:DtdivvbLptemp}), (\ref{eq:divvbW1ptd}) and (\ref{eq:qW2ptemp}), and (\ref{eq:divvbLip}) will be clear once we obtain (\ref{eq:vBH3}). $\Box$\\[2mm]
\textit{Proof of (\ref{eq:Dt2divvap})-(\ref{eq:gradDtqap}): }The proof again relies on Lemma \ref{lem:pardiv}, and we need to derive a parabolic equation for $U=D_t^2\divg v$. Taking $B=0$ in (\ref{eq:parVpm}) and dropping the subscript $\pm$, we get
\begin{equation}\label{eq:parV}
D_tV-\divg(\calt\nabla V)=D_t\divg G+\divg F+R,
\end{equation}
where $V=D_t\divg v$ and 
\[
G=2\divg v\nabla\calt-\sum_{j=1}^3\nabla_j\calt(\nabla_jv+\nabla v_j). 
\]
For $\divg F,R$, it is better to cancel $\divg(\divg v\calt\nabla\divg v)$ in $\divg F$ with $-\divg v\divg(\calt\nabla\divg v)$ in $R$ and rewrite them as
\[
F=[D_t,\calt\nabla]\divg v-\sum_{j=1}^3\calt\nabla_j\divg v\nabla_jv,
\]
\[
R=-D_t(\divg v^2)+\calt|\nabla\divg v|^2. 
\]
Applying $D_t$ to both sides in (\ref{eq:parV}) and using (\ref{eq:commdivform}), we get
\[
D_tU-\divg(\calt\nabla U)=D_t^2G+D_t\divg F+\divg H+S,
\]
where
\[
H^j=[D_t,\calt\nabla_j]V-\calt\nabla V\cdot\nabla v^j,\quad S=D_tR+\calt\nabla V\cdot\nabla\divg v. 
\]
We need more elliptic bounds: 
\begin{equation}\label{eq:DtqW1p}
\Vert\nabla D_tq\Vert_{L^p(\Omega)}\leq C(K,M,L_0,L_1)(1+\Vert U\Vert_{L^4(\Omega)}),\quad p<\infty,
\end{equation} 
\begin{equation}\label{eq:VW1p}
\Vert\nabla V\Vert_{L^p(\Omega)}\leq C(K,M,L_0,L_1)(1+\Vert U\Vert_{L^4(\Omega)}),\quad p<\infty. 
\end{equation}
To prove (\ref{eq:DtqW1p}), we use (\ref{eq:divTdq}). When $B=0$, it becomes
\[
-\divg(\calt\nabla q)=D_t\divg v+\nabla_iv^j\nabla_jv^i. 
\]
Apply $D_t$ to both sides, commute derivatives and use (\ref{eq:qW2p}), then we get
\[
\Vert D_t\divg(\calt\nabla q)\Vert_{L^4(\Omega)}\leq C(K,M,L_0,L_1)(1+\Vert U\Vert_{L^4(\Omega)}). 
\]
Next, by 
\begin{equation}\label{eq:Dtqellip}
\begin{aligned}
\divg(\calt\nabla D_tq)=\divg\{[\calt\nabla,D_t]q\}+\nabla v^j\cdot\nabla_j(\calt\nabla q)+D_t\divg(\calt\nabla q)
\end{aligned}
\end{equation}
and (\ref{eq:W1pest}), (\ref{eq:DtqW1p}) thus follows. The proof of (\ref{eq:VW1p}) is similar, but this time using (\ref{eq:parV}). As shown in the proof of (\ref{eq:DtdivvbLp}), $D_t\divg G$ has the form $\divg\xi+\eta$ with $\Vert\xi\Vert_{L^p(\Omega)},\ p<\infty$ and $\Vert\eta\Vert_{L^4(\Omega)}$ controlled by $C(K,M,L_0,L_1)$. Then (\ref{eq:VW1p}) is clear from (\ref{eq:W1pest}). 

Now we are in a position to show that
\[
D_tU-\divg(\calt\nabla U)=\divg\mathcal F+\mathcal R
\]
with 
\[
\Vert\mathcal F\Vert_{L^4(\Omega)}+\Vert\mathcal R\Vert_{L^4(\Omega)}\leq C(K,M,L_0,L_1)(1+\Vert U\Vert_{L^4(\Omega)}).
\]
The reader should not confuse the $\mathcal F,\mathcal R$ here with the ones in previous proofs. By the newly established estimates, we may well attribute $S$ to $\mathcal R$ and $H$ to $\mathcal F$. For $D_t\divg F$, we commute derivatives. The commutator $[\divg,D_t]F=\nabla_j(F\cdot\nabla v^j)-F\cdot\nabla\divg v$ can be controlled as $\divg\mathcal F+\mathcal R$, and by $D_tv=-\calt\nabla q$, (\ref{eq:VW1p}), we have control of $D_tF$ as $\mathcal F$. The toughest terms are $D_t^2\divg G$. We deduce from (\ref{eq:commdivform}) that
\begin{equation}\label{eq:comm2divform}
\begin{aligned}
[\divg,D_t^2]G=2\nabla_j(D_tG\cdot\nabla v^j)-2D_tG\cdot\nabla\divg v+\nabla_j(G\cdot D_t\nabla v^j)-G\cdot D_t\nabla\divg v\\
-\nabla_j(G\cdot\nabla v^k\nabla_kv^j)+G\cdot\nabla v^k\nabla_k\divg v. 
\end{aligned}
\end{equation}
Above arguments show that $[\divg,D_t^2]G$ can be controlled as $\divg\mathcal F+\mathcal G$ as well. After distributing $D_t^2$ onto $G$ and commuting derivatives, the most annoying term is 
\[
\nabla_i(\nabla\calt\cdot\nabla D_t^2v^i)+\divg(\nabla_j\calt D_t\nabla D_tv^j):=I_1+I_2. 
\]
For $I_1$, we write 
\[
I_1=\divg(\nabla_jD_t^2v^j\nabla\calt)+\divg(D_t^2v^j\nabla_j\nabla\calt)-\nabla_j(D_t^2v^j\divg v),
\]
noting $\Delta\calt=\divg v$. It is easy to see that $D_t^2v^j\nabla_j\nabla\calt$, $D_t^2v\divg v$ can be controlled like $\mathcal F$. For $\nabla_jD_t^2v^j$, we again commute $\nabla_j$ with $D_t^2$ and get $D_t^2\divg v$, which is just $U$. For $I_2$, write
\[
D_t\nabla D_tv^j=-D_t\nabla(\calt\nabla_jq)=-D_t\nabla_j(\calt\nabla q)+D_t(\nabla_j\calt\nabla q-\nabla\calt\nabla_jq),
\]
\[
\nabla_j\calt D_t\nabla_j(\calt\nabla q)=\nabla_j(\nabla_j\calt D_t(\calt\nabla q))-\nabla_j\nabla_j\calt D_t(\calt\nabla q)-\nabla_j\calt[\nabla_j,D_t](\calt\nabla q). 
\]
Then 
\[\begin{aligned}
I_2&=-\sum_{i,j}\nabla_i\nabla_j(\nabla_j\calt D_t(\calt\nabla_i q))+\divg\mathcal F\\
&=-\sum_j\nabla_j(\nabla_j\calt\divg D_t(\calt\nabla q))+\divg\mathcal F,
\end{aligned}\]
where $\mathcal F$ is a generic term with desired control. We have seen that $D_t\divg(\calt\nabla q)$ and $[D_t,\divg](\calt\nabla q)$ are controlled in $L^4(\Omega)$ by $C(K,M,L_0,L_1)(1+\Vert U\Vert_{L^4(\Omega)})$. Till now we establish all estimates needed for Lemma \ref{lem:pardiv} and conclude that $\Vert D_t^2\divg v\Vert_{L^4(\Omega)}\leq C(K,M,L_0,L_1)(1+\Vert D_t^2\divg v\Vert_{L^4(\Omega)}|_{t=0})$. Inequality (\ref{eq:gradDtqap}) follows from (\ref{eq:Dtqellip}), (\ref{eq:gradest}) and (\ref{eq:vBH3}). $\Box$

\subsection{Elliptic estimates}\label{sec:ellipest}

In this subsection we derive the following estimates.

\begin{proposition}\label{prop:ellipest}
(i) For $v,B$ we have
\begin{equation}\label{eq:D3vB}
\sum_{l=0}^2\sum_{k=0}^{3-l}(\Vert\nabla^kD_t^lv\Vert_{L^2(\Omega)}^2+\Vert\nabla^kD_t^lB\Vert_{L^2(\Omega)}^2)\leq C(\mathcal K)\mathcal E. 
\end{equation}
(ii) For $\divg v,\divg b$ we have
\begin{equation}\label{eq:divvbH3}
\Vert\nabla^3\divg v\Vert_{L^2(\Omega)}+\Vert\nabla^3\divg b\Vert_{L^2(\Omega)}\leq C(\mathcal K)\sqrt\mathcal E\log\mathcal E,
\end{equation}
\begin{equation}\label{eq:D3divvb}
\sum_{l=1}^2\sum_{k=0}^{3-l}(\Vert\nabla^kD_t^l\divg v\Vert_{L^2(\Omega)}^2+\Vert\nabla^kD_t^l\divg b\Vert_{L^2(\Omega)}^2)\leq C(\mathcal K)\mathcal E,
\end{equation}
and 
\begin{equation}\label{eq:VpmW2p}
\Vert\nabla^2D_-\divg(v+2b)\Vert_{L^p(\Omega)}^2+\Vert\nabla^2D_+\divg(v-2b)\Vert_{L^p(\Omega)}^2\leq C(\mathcal K)\mathcal E,\quad p<6. 
\end{equation}
(iii) For $q,N$ we have
\begin{equation}\label{eq:D3q}
\sum_{k=0}^3\Vert\nabla^kq\Vert_{L^2(\Omega)}^2+\sum_{k=0}^2\Vert\nabla^kD_tq\Vert_{L^2(\Omega)}^2\leq C(\mathcal K)\mathcal E,
\end{equation}
\begin{equation}\label{eq:D2qNbdry}
\Vert\Pi(\nabla^3q)\Vert_{L^2(\partial\Omega)}^2+\Vert\nabla^2q\Vert_{L^2(\partial\Omega)}^2+\Vert\overline\nabla^2N\Vert_{L^2(\partial\Omega)}^2\leq C(\mathcal K)\mathcal E,
\end{equation}
and
\begin{equation}\label{eq:D3qbdry}
\Vert\nabla^3q\Vert_{L^2(\partial\Omega)}+\Vert\nabla^2D_tq\Vert_{L^2(\partial\Omega)}\leq C(\mathcal K)\sqrt\mathcal E\log\mathcal E,
\end{equation}
\begin{equation}\label{eq:PiD3Dtqbdry}
\Vert\Pi(\nabla^3D_tq)\Vert_{L^2(\partial\Omega)}\leq C(\mathcal K)\sqrt\mathcal E(\log\mathcal E+\Vert\nabla_ND_tq\Vert_{L^\infty(\partial\Omega)}). 
\end{equation}
(iv) For $\calt$ we have
\begin{equation}\label{eq:TH3}
\Vert(1+|\nabla v|)\nabla^3\calt\Vert_{L^2(\Omega)}^2\leq C(\mathcal K)\mathcal E. 
\end{equation}
\end{proposition}
The proof is presented in several lemmas. 

\begin{lemma}
\begin{equation}\label{eq:vBH3}
\sum_{k=1}^3(\Vert\nabla^kv\Vert_{L^2(\Omega)}^2+\Vert\nabla^kB\Vert_{L^2(\Omega)}^2)\leq C(\mathcal K)\mathcal E,
\end{equation}
\begin{equation}\label{eq:divvbW26}
\Vert\nabla^2\divg v\Vert_{L^6(\Omega)}^2+\Vert\nabla^2\divg b\Vert_{L^6(\Omega)}^2\leq C(\mathcal K)\mathcal E.
\end{equation}
\end{lemma}
\textit{Proof: }The proof is done by iterating the pointwise estimate (\ref{eq:ptHodge}):
\[%\begin{equation}\label{eq:ptHodge}
|\nabla^ku|^2\leq C(|\nabla^{k-1}\divg u|^2+|\nabla^{k-1}\curl u|^2+\delta^{ij}Q(\nabla^ku_i,\nabla^ku_j)). 
\]%\end{equation}
Indeed, from (\ref{eq:divvbW1p}) and (\ref{eq:ptHodge}) we bound $\Vert\nabla^2(v,B)\Vert_{L^2(\Omega)}$, then using (\ref{eq:DtdivvbLp}) and (\ref{eq:W2pest}) we bound $\Vert\nabla^2\divg(v\pm2b)\Vert_{L^2(\Omega)}$. Again using (\ref{eq:ptHodge}), (\ref{eq:W2pest}) we obtain the desired estimates. $\Box$

\begin{lemma}
\[
\Vert\nabla D_t^2\divg v\Vert_{L^2(\Omega)}^2+\Vert\nabla D_t^2\divg b\Vert_{L^2(\Omega)}^2\lesssim\mathcal E_3.
\]
\[
\Vert\nabla^2 D_t\divg v\Vert_{L^2(\Omega)}^2+\Vert\nabla^2 D_t\divg b\Vert_{L^2(\Omega)}^2\lesssim\mathcal E_3.
\]
\[
\Vert\nabla^2 D_-\divg(v+2b)\Vert_{L^p(\Omega)}^2+\Vert\nabla^2 D_+\divg(v-2b)\Vert_{L^p(\Omega)}^2\lesssim\mathcal E_3,\quad p<6.
\]
\end{lemma}
\textit{Proof: }The proof relies on the identities (\ref{eq:Dtdivvrec}), (\ref{eq:Dtdivbrec}) and
\begin{equation}\label{eq:Dt2divvrec}
\begin{aligned}
D_+D_-\divg(v+2b)+D_-D_+\divg(v-2b)=2D_t^2\divg v-2b^k\nabla_k(b^j\nabla_j\divg v)\\
+2(\divg v)b^k\nabla_k\divg b,
\end{aligned}
\end{equation}
\begin{equation}\label{eq:Dt2divbrec}
\begin{aligned}
D_+D_-\divg(v+2b)-D_-D_+\divg(v-2b)=4D_t^2\divg b-4b^k\nabla_k(b^j\nabla_j\divg b)\\
+(\divg v)b^k\nabla_k\divg v.
\end{aligned}
\end{equation}
The quantities on the left come from the energy. To control the various remainder terms on the right, we make use of (\ref{eq:pardivvb}). By $E_2^{c\pm}$ we control $\Vert D_-\divg(v+2b)\Vert_{L^2(\Omega)}^2,\Vert D_+\divg(v-2b)\Vert_{L^2(\Omega)}^2$, then by (\ref{eq:pardivvb}) we have
\begin{equation}\label{eq:DeltadivvbH1}
\Vert\nabla\Delta\divg(v+2b)\Vert_{L^2(\Omega)}^2+\Vert\nabla\Delta\divg(v-2b)\Vert_{L^2(\Omega)}^2\le C(\mathcal K)\mathcal E. 
\end{equation}
Now we are able to bound $b^k\nabla_k\divg v,b^k\nabla_k\divg b$ in $H^2(\Omega)$. Indeed, they vanish on $\partial\Omega$, and
\[
\Vert\Delta(b^k\nabla_k\divg v)\Vert_{L^2(\Omega)}\leq\Vert[\Delta,b^k\nabla_k]\divg v\Vert_{L^2(\Omega)}+\Vert b^k\nabla_k\Delta\divg v\Vert_{L^2(\Omega)}\le C(\mathcal K)\sqrt{\mathcal E},
\]
\[
\Vert\Delta(b^k\nabla_k\divg b)\Vert_{L^2(\Omega)}\leq\Vert[\Delta,b^k\nabla_k]\divg b\Vert_{L^2(\Omega)}+\Vert b^k\nabla_k\Delta\divg b\Vert_{L^2(\Omega)}\le C(\mathcal K)\sqrt{\mathcal E},
\]
then from elliptic estimates we derive
\begin{equation}\label{eq:bdivvbH2}
\Vert\nabla^2(b^k\nabla_k\divg v)\Vert_{L^2(\Omega)}^2+\Vert\nabla^2(b^k\nabla_k\divg b)\Vert_{L^2(\Omega)}^2\le C(\mathcal K)\mathcal E. 
\end{equation}
Using (\ref{eq:Dt2divvrec}), (\ref{eq:Dt2divbrec}) we then obtain the desired bound for $\nabla D_t^2\divg v$, $\nabla D_t^2\divg b$. 

The second statement follows the third combined with (\ref{eq:Dtdivvrec}), (\ref{eq:Dtdivbrec}) and (\ref{eq:bdivvbH2}). To bound $D_-\divg(v+2b)$, $D_+\divg(v-2b)$ in $W^{2,p}(\Omega)$, we make use of (\ref{eq:parVpm}). Having $\Vert D_tV_\pm\Vert_{L^6(\Omega)}\leq C(\mathcal K)\sqrt\mathcal E$ in hand, we deduce from (\ref{eq:W1pest}) that
\begin{equation}\label{eq:DtdivvbW16}
\Vert\nabla D_-\divg(v+2b)\Vert_{L^6(\Omega)}+\Vert\nabla D_+\divg(v-2b)\Vert_{L^6(\Omega)}\le C(\mathcal K)\sqrt\mathcal E; 
\end{equation}
indeed, many of the related calculations have been done in the proof of (\ref{eq:DtdivvbLp}). As a consequence, we have
\begin{equation}\label{eq:DeltaqW16} 
\Vert\nabla\Delta q\Vert_{L^p(\Omega)}\le C(\mathcal K)\sqrt{\mathcal E},\quad p<6
\end{equation}
due to (\ref{eq:divTdq}). Now we are able to control the right-hand side in (\ref{eq:parVpm}). From (\ref{eq:divvbW26}), (\ref{eq:divvbW1p}) it follows that $\Vert \divg F_\pm\Vert_{L^p(\Omega)}+\Vert R_\pm\Vert_{L^p(\Omega)}\leq C(\mathcal K)\sqrt\mathcal E$, $p<6$. For $D_-\divg G_+$, we first commute $D_-,\divg$ and then note from (\ref{eq:DtdivvbW16}), (\ref{eq:Dmpvb}), (\ref{eq:DeltaqW16}) that $\Vert\divg D_-G_+\Vert_{L^p(\Omega)}+\Vert\divg D_+G_-\Vert_{L^p(\Omega)}\leq C(\mathcal K)\sqrt\mathcal E$, $p<6$. By (\ref{eq:W2pest}) the proof is then completed. $\Box$

\begin{lemma}
\[
\Vert\Pi(\nabla^3 q)\Vert_{L^2(\partial\Omega)}^2+\Vert\nabla^2 q\Vert_{L^2(\partial\Omega)}^2+\Vert \overline\nabla^2N\Vert_{L^2(\partial\Omega)}^2\leq C(\mathcal K)\mathcal E,
\]
\[
\Vert\nabla^3 q\Vert_{L^2(\Omega)}^2+\Vert(1+|\nabla v|)\nabla^3\calt\Vert_{L^2(\Omega)}^2\leq C(\mathcal K)\mathcal E,
\]
\[
\Vert\nabla^3\divg v\Vert_{L^2(\Omega)}+\Vert\nabla^3\divg b\Vert_{L^2(\Omega)}\leq C(\mathcal K)\sqrt\mathcal E\log\mathcal E. 
\]
\end{lemma}
\textit{Proof: }The first inequality follows from the definition of $E_3^a$ and (\ref{eq:proj}), (\ref{eq:bdryHr}), combined with (\ref{eq:qW2p}), (\ref{eq:DeltaqW16}). The rest follow from (\ref{eq:diri}), (\ref{eq:wHr}), (\ref{eq:tr}) combined with (\ref{eq:DeltadivvbH1}), (\ref{eq:divvbLip}). $\Box$

\begin{lemma}
\begin{equation}\label{eq:DeltaqDtq}
\Vert\nabla^2\Delta q\Vert_{L^2(\Omega)}+\Vert\nabla\Delta D_tq\Vert_{L^2(\Omega)}\leq C(\mathcal K)\sqrt\mathcal E\log\mathcal E,
\end{equation}
\[
\Vert\nabla^3 q\Vert_{L^2(\partial\Omega)}+\Vert\nabla^2 D_tq\Vert_{L^2(\partial\Omega)}\leq C(\mathcal K)\sqrt\mathcal E\log\mathcal E,
\]
\[
\Vert\Pi(\nabla^3 D_tq)\Vert_{L^2(\partial\Omega)}\leq C(\mathcal K)\sqrt\mathcal E(\log\mathcal E+\Vert\nabla_ND_tq\Vert_{L^\infty(\partial\Omega)}). 
\]
\end{lemma}
\textit{Proof: }By the reasoning in the preceding lemma, it suffices to prove the first inequality. The estimate of $\nabla^2\Delta q$ is immediate from (\ref{eq:divTdq}) combined with (\ref{eq:vBLip}), (\ref{eq:qW2p}) and the estimate of $\nabla^3(q,v,B),\ \nabla^2D_t\divg v$. For the control of $\nabla\Delta D_tq$, we again look at (\ref{eq:divTdq}). Applying $D_t$ to both sides and commuting derivatives, we get
\[\begin{aligned}
-\calt\Delta D_tq=D_t^2\divg v+[D_t,\calt\Delta]q&+D_t(\nabla q\cdot\nabla\calt)\\
&+D_t\{\nabla_iv^j\nabla_jv^i-\nabla_i(\calt B^j)\nabla_jB^i\}. 
\end{aligned}\]
Then it is clear that $\Vert \Delta D_tq\Vert_{L^p(\Omega)}\leq C(\mathcal K)\sqrt\mathcal E$, $p<6$,
and as a consequence, 
\begin{equation}\label{eq:DtqW2p}
\Vert\nabla^2D_tq\Vert_{L^p(\Omega)}\leq C(\mathcal K)\sqrt\mathcal E,\quad p<6. 
\end{equation}
Further applying $\nabla$ to both sides, by a lengthy calculation we obtain the desired estimate. $\Box$\\[2mm]
\textit{Proof of (\ref{eq:D3vB})-(\ref{eq:TH3}): }Most of the proofs are included in the preceding lemmas. The rest work is using (\ref{eq:LM}) for converting derivatives. Details are omitted. $\Box$

\section*{Acknowledgements}
This work is partially supported by an NSFC Grant 12171267 and a New Cornerstone Investigator Program 100001127 through my Ph. D. supervisor, Professor Huihui Zeng. The author would like to thank Huihui Zeng for suggesting this problem and for many helpful discussions.

\appendix

\section{Sobolev inequalities}

Let $\Omega\subset\mathbb R^n$ be a bounded domain with $|\theta|+\frac{1}{\iota_0}\leq K$. 

\begin{lemma}\label{lem:soboineq}
(i) 
\[
\Vert u\Vert_{L^{np/(n-kp)}(\Omega)}\leq C(K)\sum_{j=0}^k\Vert\nabla^ju\Vert_{L^p(\Omega)},\quad 1\leq p<\frac{n}{k},
\]\[
\Vert u\Vert_{L^\infty(\Omega)}\leq C(K)\sum_{j=0}^k\Vert\nabla^ju\Vert_{L^p(\Omega)},\quad p>\frac{n}{k}.
\]
(ii) Let either $D^j=\nabla^j,\ j=1,...,k$ or $D^j=\overline\nabla^j,\ j=1,...,k$. Then we have
\[
\Vert u\Vert_{L^{(n-1)p/(n-1-kp)}(\partial\Omega)}\leq C(K)\sum_{j=0}^k\Vert D^ju\Vert_{L^p(\partial\Omega)},\quad 1\leq p<\frac{n-1}{k},
\]\[
\Vert u\Vert_{L^\infty(\partial\Omega)}\leq C(K)\sum_{j=0}^k\Vert D^ju\Vert_{L^p(\partial\Omega)},\quad p>\frac{n-1}{k}.
\]
\end{lemma}
\textit{Proof: }These inequalities with $D^j=\nabla^j$ can be deduced from Lemmas 3.6, A.2-A.4 in \cite{CLin}. The arguments therein also apply for $D^j=\overline\nabla^j$. $\Box$

\begin{corollary}\label{cor:prod}
(i) Let $r_1+\cdots+r_m\geq\frac{n}{2}(m-1)$. Then we have
\begin{equation}\label{eq:prodint}
\Vert\prod_{i=1}^mf_i\Vert_{L^2(\Omega)}\leq C(K)\prod_{i=1}^m\left(\sum_{j=0}^{r_i}\Vert\nabla^jf_i\Vert_{L^2(\Omega)}+\sigma\sigma_i\Vert f_i\Vert_{L^\infty(\Omega)}\right),
\end{equation}
where
\begin{equation}\label{eq:sigmaint}
\sigma=\begin{cases}1&\sum r_i=\frac{n}{2}(m-1),\\0&\mathrm{otherwise},\end{cases}\qquad\sigma_i=\begin{cases}1& r_i=\frac{n}{2},\\0&\mathrm{otherwise}.\end{cases}
\end{equation}
(ii) Let $r_1+\cdots+r_m\geq\frac{n-1}{2}(m-1)$ and let $D^{\leq r_i}$ be either $\nabla^{\leq r_i}$ or $\overline\nabla^{\leq r_i}$. Then we have
\begin{equation}\label{eq:prodbdry}
\Vert\prod_{i=1}^mf_i\Vert_{L^2(\partial\Omega)}\leq C(K)\prod_{i=1}^m\left(\sum_{j=0}^{r_i}\Vert D^jf_i\Vert_{L^2(\partial\Omega)}+\sigma\sigma_i\Vert f_i\Vert_{L^\infty(\partial\Omega)}\right),
\end{equation}
where
\begin{equation}\label{eq:sigmabdry}
\sigma=\begin{cases}1&\sum r_i=\frac{n-1}{2}(m-1),\\0&\mathrm{otherwise},\end{cases}\qquad\sigma_i=\begin{cases}1& r_i=\frac{n-1}{2},\\0&\mathrm{otherwise}.\end{cases}
\end{equation}
\end{corollary}
\textit{Proof: }(i) The condition implies that there are $p_i\in[1,\infty]$ such that
\[
\frac{1}{p_1}+\cdots+\frac{1}{p_m}=\frac{1}{2},\qquad \frac{1}{p_i}\geq\frac{1}{2}-\frac{r_i}{n}
\]
such that $\sigma\sigma_i=1$ if $r_i=\frac{n}{2}$ and $p_i=\infty$. Then (\ref{eq:prodint}) follows Lemma \ref{lem:soboineq} and H\"older inequality. \\
(ii) The proof is similar. $\Box$

\begin{corollary}\label{cor:nonlin}
Let $U\subset\mathbb R^N$ be a bounded region and $\phi:U\to\mathbb R$ be a smooth function. 
Let $f:[0,T]\times\overline\Omega\to U$ and $\delta=d(f([0,T]\times\overline\Omega),\mathbb R^N\backslash U)$. With the notation in (\ref{eq:rsnorm}) we have
\begin{equation}\label{eq:nonlin1}
\Vert\phi(f)\Vert_{r,s}\leq C(K,\phi,\delta^{-1},\Vert f\Vert_{r,s}),\quad r>\frac{n}{2},
\end{equation}
and
\begin{equation}\label{eq:nonlin2}
\Vert\phi(f)\Vert_{r,s}\leq C(K,\phi,\delta^{-1},\Vert f\Vert_{r-1,s-1})\Vert f\Vert_{r,s},\quad r>\frac{n}{2}+1. 
\end{equation}
\end{corollary}
\textit{Proof: }Let $\Gamma^r$ be a composition of $r$ elements of $\{D_t,\nabla\}$ where $D_t$ appears at most $s$ times. By chain rule, $\Gamma^r\phi(f)$ is a sum of the following general terms:
\[
C(f)\Gamma^{r_1}f\cdots\Gamma^{r_m}f,
\]
where $C(f)$ is a smooth function in $f$, $\sum r_i=r$ and $D_t$ appears at most $s$ times in the whole expression. Since $\sum_{i=1}^m(r-r_i)>\frac{n}{2}(m-1)$ if $r>\frac{n}{2}$ and $\sum_{i\neq j}(r-1-r_i)+(r-r_j)>\frac{n}{2}(m-1)$ if $r>\frac{n}{2}+1$, where $j$ is such that $\Gamma^{r_j}$ contains most $D_t$'s, the results simply follow from (\ref{eq:prodint}). $\Box$

\section{Elliptic and interpolation estimates}

Let $\Omega\subset\mathbb R^3$ be a bounded $C^2$ domain with $|\theta|+\frac{1}{\iota_0}\leq\frac{K}{10}$. The following estimate for is classical; see \cite{MulInt}. 

\begin{theorem}\label{thm:diriest}
Let $u$ solve the problem
\[
\left\{\begin{aligned}
\Delta u&=\divg w+f\quad&\text{in }\Omega,\\
u&=0\quad&\text{on }\partial\Omega. 
\end{aligned}\right.
\]
Then for $1<p<\infty$ and $q>3$ we have
\begin{equation}\label{eq:W1pest}
\Vert\nabla u\Vert_{L^p(\Omega)}\leq C(K)\Vert w\Vert_{L^p(\Omega)}+\Vert f\Vert_{L^q(\Omega)}
\end{equation}
and
\begin{equation}\label{eq:W2pest}
\Vert\nabla^2 u\Vert_{L^p(\Omega)}\leq C(K)\Vert\divg w+f\Vert_{L^p(\Omega)}. 
\end{equation}
\end{theorem}
The Green function for the Dirichlet problem enjoys the following estimate; see \cite{GWid}. 

\begin{theorem}\label{thm:grfunc}
There exists a function $G(x,y)$ such that
\[
\phi(x)=\int_{\Omega}G(x,y)(-\Delta \phi(y))dy=\int_{\Omega}\nabla_yG(x,y)\nabla \phi(y)dy,\quad x\in\Omega
\]
for $\phi\in H_0^1\cap H^2(\Omega)$. Moreover, $G(x,\cdot)=0$ on $\partial\Omega$, $G(x,y)=G(y,x)$ and
\begin{equation}\label{eq:grest}
0\leq G(x,y)\leq \frac{1}{4\pi|x-y|},\quad |\nabla_xG(x,y)|\leq \frac{C(K)}{|x-y|^2},\quad |\nabla_y\nabla_xG(x,y)|\leq \frac{C(K)}{|x-y|^3}.
\end{equation}
\end{theorem}
A corollary to the above estimate is (see \cite{Gins}): 

\begin{theorem}\label{thm:gradest}
Let $u$ solve the problem
\[
\left\{\begin{aligned}
\Delta u&=\divg w+f\quad&\text{in }\Omega,\\
u&=0\quad&\text{on }\partial\Omega. 
\end{aligned}\right.
\]
Then it holds that
\begin{equation}\label{eq:gradest}
\Vert\nabla u\Vert_{L^\infty(\Omega)}\leq C(K)\left[\Vert f\Vert_{L^4(\Omega)}+\Vert w\Vert_{L^\infty(\Omega)}\bigg(1+\log^+\frac{\Vert w\Vert_{H^2(\Omega)}}{\Vert w\Vert_{L^\infty(\Omega)}}\bigg)\right].
\end{equation}
\end{theorem}
Another corollary is

\begin{lemma}\label{lem:gradest2}
Let $u$ be harmonic in $\Omega$ and $\alpha\in(0,1]$. Then we have
\begin{equation}\label{eq:gradLp}
\Vert\nabla u\Vert_{L^p(\Omega)}\leq C(K,\vol\Omega)\Vert\overline\nabla u\Vert_{L^\infty(\partial\Omega)}
\end{equation}
and
\begin{equation}\label{eq:gradest2}
\Vert\nabla u\Vert_{L^\infty(\Omega)}\leq C(K,\vol\Omega)\Vert\overline\nabla u\Vert_{L^\infty(\partial\Omega)}\Big(1+\log^+\frac{\Vert\overline{\nabla} u\Vert_{C^\alpha(\partial\Omega)}}{\Vert\overline\nabla u\Vert_{L^\infty(\partial\Omega)}}\Big). 
\end{equation}
\end{lemma}
\textit{Proof: }We make use of the Green representation formula: 
\begin{equation}\label{eq:grrep}
\nabla u(x)=\int_{\partial\Omega}\nabla_x\partial_{N_y}G(x,y)u(y)dS(y). 
\end{equation}
Since $\int_{\partial\Omega}\partial_{N_y}G(x,y)dS(y)\equiv1$, we also have
\[
\nabla u(x)=\int_{\partial\Omega}\nabla_x\partial_{N_y}G(x,y)(u(y)-u(\bar x))dS(y),
\]
where $\bar x\in\partial\Omega$ is such that $|x-\bar x|=d(x,\partial\Omega)$. Notice that
\begin{equation}\label{eq:2ptest}
|u(y)-u(\bar x)|\leq \Vert\overline\nabla u\Vert_{L^\infty(\partial\Omega)}d_\zeta(y,\bar x), 
\end{equation}
where $d_\zeta$ is the geodesic distance on $\partial\Omega$. This allows us to control $\text{osc}_{\partial\Omega}u$. Indeed, by the 3-time covering lemma and $\iota_0\geq \frac{10}{K}$, there exist $m$ balls of radius $\frac{1}{3K}$ covering $\partial\Omega$ with $\frac{m}{K^2}\leq C\text{Area}\,\partial\Omega\leq CK\vol\Omega$. Hence, every pair of points on $\partial\Omega$ can be joint by a curve of length less or equal than $CK^2\vol\Omega$. Then by (\ref{eq:2ptest}) we have
\begin{equation}\label{eq:oscest}
\sup_{y,z\in\partial\Omega}|u(y)-u(z)|\leq C(K,\vol\Omega)\Vert\overline\nabla u\Vert_{L^\infty(\partial\Omega)}. 
\end{equation}
Noting that $|y-x|\geq C^{-1}(|y-\bar x|+d(x,\partial\Omega))$, $|y-\bar x|\leq\frac{1}{K}$ and using (\ref{eq:grest}) and the geodesic polar coordinates, we have
\[\begin{aligned}
|\nabla u(x)|\leq&\left(\int_{\partial\Omega\cap B_\frac{1}{K}(\bar x)}+\int_{\partial\Omega\backslash B_\frac{1}{K}(\bar x)}\right)|\nabla_x\partial_{N_y}G(x,y)(u(y)-u(\bar x))|dS(y)\\
\leq&\int_0^\frac{1}{K}C(K)\frac{\Vert\overline\nabla u\Vert_{L^\infty(\partial\Omega)}r}{(r+d(x,\partial\Omega))^3}rdr+C(K,\vol\Omega)\Vert\overline\nabla u\Vert_{L^\infty(\partial\Omega)}\\
\leq&C(K,\vol\Omega)\Vert\overline\nabla u\Vert_{L^\infty(\partial\Omega)}\left(1+\log^+\frac{1}{d(x,\partial\Omega)}\right). 
\end{aligned}\]
By considering the foliation consisting of level sets of $d(x,\partial\Omega)$, we obtain (\ref{eq:gradLp}). 

To prove (\ref{eq:gradest2}), we need to subtract the linear part from (\ref{eq:grrep}). Let $x\in\Omega$ be fixed and $a=\overline\nabla u(\bar x)$, $w(y)=u(y)-u(\bar x)-a\cdot y$. Since linear functions are harmonic, we have $\int_{\partial\Omega}\partial_{N_y}G(z,y)(u(\bar x)+a\cdot y)dS(y)=u(\bar x)+a\cdot z$, $z\in\Omega$, hence
\[
\nabla u(x)=\int_{\partial\Omega}\nabla_x\partial_{N_y}G(x,y)w(y)dS(y)+O(\Vert\overline\nabla u\Vert_{L^\infty(\Omega)}). 
\]
Next, for $r\in(0,\frac{1}{K})$, decompose the integral into three parts: 
\[
\bigg(\int_{\partial\Omega\cap B_r(\bar x)}+\int_{\partial\Omega\cap B_{\frac{2}{K}}\backslash B_r(\bar x)}+\int_{\partial\Omega\backslash B_{\frac{2}{K}}(\bar x)}\bigg)\nabla_x\partial_{N_y}G(x,y)w(y)dS(y):=I_1+I_2+I_3. 
\]
Denote $\delta=d(x,\partial\Omega)$ and note that $|w(y)|\leq|y-\bar x|^{1+\alpha}\Vert\overline u\Vert_{C^\alpha(\partial\Omega)}$, $|y-\bar x|\leq\frac{1}{K}$. By (\ref{eq:grest}) and using the geodesic polar coordinates, we have
\[\begin{aligned}
|I_1|\leq& C(K)\Vert\overline\nabla u\Vert_{C^{\alpha}(\partial\Omega)}\int_{\partial\Omega\cap B_r(\bar x)}\frac{1}{|x-y|^3}|y-\bar x|^{1+\alpha}dS(y)\\
\leq& C(K)\Vert\overline\nabla u\Vert_{C^{\alpha}(\partial\Omega)}\int_0^r\frac{t^{1+\alpha}}{(\delta+t)^3}tdt\\
\leq& C(K)\Vert\overline\nabla u\Vert_{C^{\alpha}(\partial\Omega)}r^\alpha,
\end{aligned}\]
\[\begin{aligned}
|I_2|\leq& C(K)\Vert\overline\nabla u\Vert_{L^\infty(\partial\Omega)}\int_{\partial\Omega\cap B_{\frac{2}{K}}\backslash B_r(\bar x)}\frac{1}{|x-y|^3}|y-\bar x|dS(y)\\
\leq& C(K)\Vert\overline\nabla u\Vert_{L^\infty(\partial\Omega)}\int_r^{\frac{2}{K}}\frac{t}{(\delta+t)^3}tdt\\
\leq& C(K)\Vert\overline\nabla u\Vert_{L^\infty(\partial\Omega)}\log \frac{1}{r+\delta},
\end{aligned}\]
\[\begin{aligned}
|I_3|\leq& C(K)\int_{\partial\Omega\backslash B_{\frac{2}{K}}(\bar x)}\frac{1}{|x-y|^3}|w(y)|dS(y)\\\leq& C(K,\vol\Omega)\Vert w\Vert_{L^\infty(\partial\Omega)}. 
\end{aligned}\]
As is shown before, every pair of points on $\partial\Omega$ can be joint by a curve of length less or equal than $CK^2\vol\Omega$, thus $|w(y)|\leq C(K,\vol\Omega)$. Finally, selecting $r^\alpha=\min\{\frac{1}{K^\alpha},\frac{\Vert\overline\nabla u\Vert_{L^\infty(\partial\Omega)}}{\Vert\overline\nabla u\Vert_{C^{\alpha}(\partial\Omega)}}\}$ we complete the proof. $\Box$

\begin{lemma}
(i) For $p,q,r\in[2,\infty]$ with $\frac{2}{p}=\frac{1}{q}+\frac{1}{r}$, it holds
\begin{equation}\label{eq:gradinterp}
\Vert\nabla f\Vert_{L^p(\Omega)}^2\leq C(K)\Vert f\Vert_{L^q(\Omega)}\Vert f\Vert_{W^{2,r}(\Omega)}. 
\end{equation}
(ii)
For $2\leq p<\infty$, it holds
\begin{equation}\label{eq:trinterp}
\Vert f\Vert_{L^p(\partial\Omega)}\leq CK^{\frac{1}{p}}\Vert f\Vert_{L^{2p-2}(\Omega)}^{1-\frac{1}{p}}\Vert f\Vert_{H^1(\Omega)}^{\frac{1}{p}}. 
\end{equation}
\end{lemma}
\textit{Proof: }(i) See Lemma A.1 in \cite{CLin}. \\
(ii) By the divergence formula, we have
\[\begin{aligned}
\Vert f\Vert_{L^p(\partial\Omega)}^p=\int_{\partial\Omega}N^iN_i|f|^pdS=&\int_\Omega\partial_i(|f|^pN^i)dx\\
=&\int_\Omega(\divg N)|f|^p+p|f|^{p-2}f\partial_Nfdx\\
\leq&CK\Vert f\Vert_{L^{2p-2}(\Omega)}^{1-\frac{1}{p}}\Vert f\Vert_{L^2(\Omega)}^{\frac{1}{p}}+p\Vert f\Vert_{L^{2p-2}(\Omega)}^{1-\frac{1}{p}}\Vert\nabla f\Vert_{L^2(\Omega)}^{\frac{1}{p}},
\end{aligned}\]
as desired. $\Box$

\section{Vector identities}\label{sec:vectid}

\begin{definition}
Let
\[
\epsilon_{ijk}=\begin{cases}1&(i,j,k)\text{ is an even permutation of }(1,2,3),\\
-1&(i,j,k)\text{ is an odd permutation of }(1,2,3),\\0&\text{otherwise}.\end{cases}
\]
Let $F^i,G^j$ be two vector fields on $\Omega\subset\mathbb R^3$. Define $(F\times G)_k=\epsilon_{ijk}F^iG^j$ and $(\curl G)_k=(\nabla\times G)_k=\epsilon^i_{\ jk}\nabla_iG^j$. More explicitly, 
\[
F\times G=(F^2G^3-F^3G^2,\ F^3G^1-F^1G^3,\ F^1G^2-F^2G^1)^T,
\]
\[
\curl G=(\nabla_2G^3-\nabla_3G^2,\ \nabla_3G^1-\nabla_1G^3,\ \nabla_1G^2-\nabla_2G^1)^T. 
\]
\end{definition}

\begin{lemma}
Let $F,G,H$ be vector fields on $\Omega\subset\mathbb R^3$. Then we have
\begin{equation}\label{eq:Deltadivcurl}
\Delta F=\nabla\divg F-\curl\curl F,
\end{equation}
\begin{equation}\label{eq:divtimes}
\divg(F\times G)=\curl F\cdot G-F\cdot\curl G,
\end{equation}
\begin{equation}\label{eq:curltimes}
\curl(F\times G)=F\divg G-G\divg F+(G\cdot\nabla)F-(F\cdot\nabla)G,
\end{equation}
\begin{equation}\label{eq:mixprod}
F\times G\cdot H=-F\times H\cdot G=-G\times F\cdot H. 
\end{equation}
\end{lemma}

\section*{Declaration}
Data availability is not applicable to this article as no new data were created or analyzed in this study. There are no conflicts of interest in relation to this manuscript.

\end{document}